\documentclass[a4paper,10pt]{article}
\usepackage{amsmath}
\usepackage{amsfonts}
\usepackage{amssymb}
\usepackage{array,longtable}
\usepackage{amscd}
\usepackage[cp1251]{inputenc}
\usepackage[english]{babel}
\usepackage[usenames]{color}
\usepackage[matrix,arrow,curve]{xy}

\begin{document}
\newtheorem{definition}{Definition}
\newtheorem{theorem}{Theorem}
\newtheorem{lemma}{Lemma}
\newtheorem{proposition}{Proposition}
\newtheorem{observation}{Observation}
\newtheorem{corollary}{Corollary}
\newcommand{\diag}{{\rm diag}}
\newcommand{\tr}{{\rm Tr}}

\begin{center}
{\large {\bf HOW TO GENERALIZE $D$-STABILITY}} \\[0.5cm]
{\bf Olga Y. Kushel} \\[0.5cm]
Shanghai University, \\ Department of Mathematics, \\ Shangda Road 99, \\ 200444 Shanghai, China \\
kushel@mail.ru
\end{center}

\begin{abstract}
In this paper, we introduce the following concept which generalizes known definitions of multiplicative and additive $D$-stability, Schur $D$-stability, $H$-stability, $D$-hyperbolicity and many others. Given a subset ${\mathfrak D} \subset {\mathbb C}$, a matrix class ${\mathcal G} \subset {\mathcal M}^{n \times n}$ and a binary operation $\circ$ on ${\mathcal M}^{n \times n}$, an $n \times n$ matrix $\mathbf A$ is called $({\mathfrak D}, {\mathcal G}, \circ)$-stable if $\sigma({\mathbf G}\circ {\mathbf A}) \subset {\mathfrak D}$ for any ${\mathbf G} \in {\mathcal G}$. Such an approach allows us to unite several well-known matrix problems and to consider common ways of their analysis. Here, we make a survey of existing results and open problems on different types of stability, study basic properties of $({\mathfrak D}, {\mathcal G}, \circ)$-stable matrices and relations between different $({\mathfrak D}, {\mathcal G}, \circ)$-stability classes.

Hurwitz stability, Schur stability, $D$-stability, $H$-stability, diagonal stability, eigenvalue clustering, Lyapunov equation.

MSC[2010] 15A48, 15A18, 15A75
\end{abstract}

\section*{Introduction} Dynamical systems theory gave rise to many classical matrix problems, basing on Lyapunov's idea that stability of a system of ODE can be established through properties of the corresponding matrix. In this connection, a number of matrix classes related to different stability types were introduced and studied, but the further development of systems analysis, including robustness and control, leads to more and more new matrix problems. Rapid growth in this area in the last decades caused some separation between classical matrix stability theory and its systems applications. Nowadays we have the following lines of research.
\begin{enumerate}
\item[-] Matrix line (for the review see \cite{HER1}, \cite{LOG}, for the examples of open stability problems see \cite{HERK2}, \cite{HOLS}, \cite{BIR}). This line mainly focuses on long-standing open problems, connected to the classes of structured matrices, introduced in 1950th--1970th and classical methods of matrix analysis.
\item[-] Systems and control line. Here, we most refer to the special monograph by Kaszkurevicz and Bhaya \cite{KAB} which collects the results on different types of diagonal stability and their applications. This research line also includes modern books in robust control theory (see \cite{BH}, \cite{BAR}, \cite{BHK2}, \cite{RUGH}), a number of books and papers describing different approaches to stability and related problems (see \cite{BGF}, \cite{GUT2}, \cite{HENG}, \cite{JONCK}, \cite{JU2}, \cite{KOG}, \cite{KUSHN}, \cite{MART}, \cite{WON}), including modern and classical polynomial methods (see \cite{JEM}, \cite{MAR}, \cite{RAS}). This highly applicable research states more general problems connected to different stability generalizations, mostly in terms of ODE systems and polynomials, and outlines some ways of their treatment.
\end{enumerate}

In this paper, we refresh the link between these two approaches, generalizing the results obtained for concrete matrix classes,  stating more general matrix problems and analyzing the ways of their solutions.

Let ${\mathcal M}^{n \times n}$ denote the set of all real $n \times n$ matrices, $\mathbf A$ be a matrix from ${\mathcal M}^{n \times n}$, and $\sigma({\mathbf A})$ denote the spectrum of $\mathbf A$ (i.e. the set of all eigenvalues of $\mathbf A$ defined as zeroes of its characteristic polynomial $f_{\mathbf A}(\lambda):= \det(\lambda{\mathbf I - {\mathbf A}})$). Here, a {\it stability region} is a set ${\mathfrak D} \subseteq {\mathbb C}$ with the property: $\lambda \in {\mathfrak D}$ implies its complex conjugate $\overline{\lambda} \in {\mathfrak D}$ (this property is needed since we study matrices with real entries). We do not assume any restrictions, (eg. connectivity or convexity) on ${\mathfrak D}$. Given a stability region ${\mathfrak D}$, an $n \times n$ matrix ${\mathbf A}$ is called {\it stable with respect to $\mathfrak D$} or simply {\it $\mathfrak D$-stable} if $\sigma({\mathbf A}) \subset {\mathfrak D}$.

The following matrix classes can be considered as examples of $\mathfrak D$-stable matrices for different stability regions.

\begin{enumerate}
\item[\rm 1.] An $n \times n$ real matrix $\mathbf A$ is called {\it Hurwitz stable} or just {\it stable} if all its eigenvalues have negative real parts (see, for example, \cite{BELL}, \cite{KAB}, \cite{MAM}). In matrix literature, positive stability is often referred: matrix is called {\it positive stable} or just {\it stable} if all its eigenvalues have positive real parts (see \cite{HER1}, \cite{JOHN1}).
\item[\rm 2.] An $n \times n$ real matrix $\mathbf A$ is called {\it Schur stable} if all its eigenvalues lie inside the unit circle, i.e. the spectral radius $\rho(\mathbf A) < 1$ (see \cite{BHK}, \cite{KAB}). This property is mostly referred as {\it convergence of matrices} (see \cite{HOJ}, p. 137).
\item[\rm 3.] An $n \times n$ real matrix $\mathbf A$ is called {\it aperiodic} if all its eigenvalues are real (see, for example, \cite{GUJU}, p. 860 and \cite{GUT2}, p. 92).
  \end{enumerate}

In the last decades, matrix and polynomial $\mathfrak D$-stability, also known as matrix (respectively, polynomial) {\it root clustering} has become an attractive area for researchers. Giving a brief overview, the theory has been developed from the simplest and the most used partial cases to more and more sophisticated and general stability regions (see, for example, \cite{GUT2}, \cite{JU2}, \cite{JU3}, \cite{JUA}, \cite{ASCM}, \cite{BAP} -\cite{BHPM}, \cite{CHIGA}, \cite{CGA}). For the classic examples of the stability regions ${\mathfrak D}$, this concept goes back to Descartes \cite{DEC} and is developed in the papers by Caushy \cite{CAU}, Sturm \cite{STU}, Hermite \cite{HERMI}, Routh \cite{ROU1}-\cite{ROU4}, Hurwitz \cite{HUR} (${\mathfrak D}$ is the open left-hand side of the complex plane) and Schur \cite{SHU2}, Cohn \cite{COH} (${\mathfrak D}$ is the open unit disk). Studying classes of matrices with real (positive) spectra (respectively, polynomials all whose zeroes are real (positive)) (for the beginning, see \cite{WAR}) can be also considered a partial case of $\mathfrak D$-stability (${\mathfrak D}$ is the real line or its positive direction).

Now we provide the central definition of this paper.

{\bf General definition}. Given a stability region ${\mathfrak D} \subset {\mathbb C}$, a matrix class ${\mathcal G} \subset {\mathcal M}^{n \times n}$ and a binary operation $\circ$ on ${\mathcal M}^{n \times n}$, we call an $n \times n$ matrix $\mathbf A$ {\it left (right) $({\mathfrak D},{\mathcal G},\circ)$-stable} if $\sigma({\mathbf G}\circ{\mathbf A}) \subset {\mathfrak D}$ (respectively, $\sigma({\mathbf A}\circ{\mathbf G}) \subset {\mathfrak D}$) for any matrix $\mathbf G$ from the class ${\mathcal G}$. Later on, we use the term "$({\mathfrak D},{\mathcal G},\circ)$-stable" for left $({\mathfrak D},{\mathcal G},\circ)$-stability.

Basic examples of stability regions are those used in the $\mathfrak D$-stability definitions, given above.
\begin{enumerate}
\item[\rm 1.] Open left-hand  (right-hand) side of the complex plane.
\item[\rm 2.] Open unit circle.
\item[\rm 3.] Real axes or its positive (negative) direction.
\end{enumerate}

To provide the examples of the class ${\mathcal G}$ of real $n \times n$ matrices, we focus on the following most important and most studied cases.
\begin{enumerate}
\item[\rm 1.] The class $\mathcal S$ of symmetric matrices, and its subclasses. Here, we pay special attention to the class $\mathcal H$ of symmetric positive definite matrices. Recall, that an $n \times n$ matrix $\mathbf A$ is called {\it positive definite (positive semidefinite)} if $\langle x, {\mathbf A}x\rangle >0$ for all nonzero $x \in {\mathbb R}^{n}$ (respectively, $\langle x, {\mathbf A}x\rangle \geq 0$ for all $x \in {\mathbb R}^{n}$) (for the definition and properties see, for example, \cite{BHAT2}).
\item[\rm 2.] The class $\mathcal D$ of diagonal matrices and its subclasses.
\item[\rm 3.] Matrices of finite rank $k$ $(1 \leq k \leq n)$.
\end{enumerate}

The most common examples of binary operations on ${\mathcal M}^{n \times n}$ are as follows.
\begin{enumerate}
\item[\rm 1.] The operation $"+"$ of matrix addition.
\item[\rm 2.] The operation $"\cdot"$ of matrix multiplication.
\item[\rm 3.] The operation $"\circ"$ of Hadamard (entry-wise) matrix multiplication (for the definitions and properties see \cite{JOHN3}).
\end{enumerate}

This paper offers a detailed account of various connections between the classes of $({\mathfrak D},{\mathcal G},\circ)$-stable matrices and gives an overview of the corresponding results. Due to the broad range of the partial cases of general $({\mathfrak D},{\mathcal G},\circ)$-stability definition which have appeared recently and their applications to the linear system theory, our main goal is to develop some common methods of studying such generalizations. In Section 1, we collect the examples of matrix classes which appeared recently in the literature, describing the corresponding stability regions $\mathfrak D$, matrix classes $\mathcal G$ and binary operations $\circ$. Then, we focus on different stability regions and their properties, focusing on so-called "Lyapunov regions" described by generalizations of the Lyapunov theorem. We consider the most used matrix classes $\mathcal G$, their basic properties and inclusion relations between them as well as the properties of binary operations on ${\mathcal M}^{n \times n}$. In Section 2, we state and prove the basic results on $({\mathfrak D},{\mathcal G},\circ)$-stability, i.e. inclusion relations among different classes and the properties of elementary operations on $({\mathfrak D},{\mathcal G},\circ)$-stable matrices. We also define a generalization of the widely used concept of Volterra--Lyapunov stability and study its relations to $({\mathfrak D},{\mathcal G},\circ)$-stability. Section 3 provides the ways of the further development of $({\mathfrak D},{\mathcal G},\circ)$-stability theory, outlines some open problems and analyzes the applications.

\section{$({\mathfrak D},{\mathcal G},\circ)$-stability: from known to new matrix classes}
\subsection{Examples of known $({\mathfrak D},{\mathcal G},\circ)$-stability classes}
 The definition of $({\mathfrak D},{\mathcal G},\circ)$-stability includes the following known matrix classes.

\begin{enumerate}
\item[\rm 1.] An $n \times n$ real matrix $\mathbf A$ is called (multiplicative) {\it $D$-stable} if ${\mathbf D}{\mathbf A}$ is stable for every positive diagonal matrix $\mathbf D$ (i.e., an $n \times n$ matrix $\mathbf D$ with positive entries on its principal diagonal, while the rest are zero). Later on, using the term "$D$-stability", we mostly refer to multiplicative $D$-stability. Here, the stability region ${\mathfrak D}$ is the open right-hand side of the complex plane, ${\mathcal G}$ is the class of positive diagonal matrices, and $\circ$ is (left-side) matrix multiplication. This matrix class was introduced in \cite{AM} in connection with some problems of mathematical economics. The literature on multiplicative $D$-stability is particularly rich due to a lot of applications (see, for example, \cite{LOG}, \cite{QR}, \cite{KAB} and references therein).
\item[\rm 2.] An $n \times n$ real matrix $\mathbf A$ is called (multiplicative) {\it $H$-stable} if ${\mathbf H}{\mathbf A}$ is stable for every symmetric positive definite matrix $\mathbf H$. Here, the stability region ${\mathfrak D}$ and the operation $\circ$ are the same as for $D$-stability, ${\mathcal G}$ is the class of symmetric positive definite matrices. This matrix class also arises in \cite{AM} under the name of $S$-stability and later studied (see \cite{CARL3}, \cite{CAS}, \cite{OSS}) under the name of $H$-stability.
\item[\rm 3.] An $n \times n$ real matrix $\mathbf A$ is called (multiplicative) {\it $D$-positive} if all the eigenvalues of ${\mathbf D}{\mathbf A}$ are positive for every positive diagonal matrix $\mathbf D$. This definition was given in \cite{BAO} in connection with stability of mechanic systems. The stability region ${\mathfrak D}$ in this case is the positive direction of the real axes, ${\mathcal G}$ is the class of positive diagonal matrices and $\circ$ is matrix multiplication, as before. It is natural to define also the following matrix class: an $n \times n$ real matrix $\mathbf A$ is called (multiplicative) {\it $D$-aperiodic} if all the eigenvalues of ${\mathbf D}{\mathbf A}$ are real for every diagonal matrix $\mathbf D$. Here, we extend the stability region ${\mathfrak D}$ to the whole real axes and also extend the class ${\mathcal G}$ to the whole class of diagonal matrices from ${\mathcal M}^{n \times n}$.
\item[\rm 4.] An $n \times n$ real matrix $\mathbf A$ is called {\it Schur $D$-stable} if ${\mathbf D}{\mathbf A}$ is Schur stable for every diagonal matrix $\mathbf D$ with $\|{\mathbf D}\|< 1$ (i.e. an $n \times n$ diagonal matrix $\mathbf D$ with $|d_{ii}|< 1$, $i = 1, \ \ldots, \ n$). Here ${\mathfrak D} = \{z \in {\mathbb C}: |z|< 1\}$, ${\mathcal G}$ is the class of diagonal matrices with $\|{\mathbf D}\|< 1$, and $\circ$ is, as before, matrix multiplication. This matrix class was defined in \cite{BHK} in connection with the study of discrete-time systems, and studied in \cite{CA2} as {\it convergent multiples}.
\item[\rm 5.] A matrix $\mathbf A$ is called {\it vertex stable} if $\rho({\mathbf D}{\mathbf A}) < 1$ for any real diagonal matrix $\mathbf D$ with $|{\mathbf D}|=1 $, i.e. with $d_{ii} = \pm1, \ i = 1, \ \ldots, \ n$ (note that this matrix class contains only a finite number of matrices). This matrix class was also defined in \cite{BHK}, in connection with Schur $D$-stability. The paper \cite{CA2} provides different proofs of the basic results on Schur $D$-stability and vertex stability. Some generalization of Schur and vertex stability is studied in \cite{GEOH}: given a matrix ${\mathbf A} \in {\mathcal M}^{n \times n}$, and a convex polyhedron ${\mathcal B} \subset {\mathcal M}^{n \times n}$, defined as the convex hull of the finite number of its {\it extreme} matrices ${\mathbf B}_1, \ \ldots, \ {\mathbf B}_N$. The stability of matrix set $${\mathcal A}: = \{{\mathbf A}{\mathbf B}: {\mathbf B} \in {\mathcal B}\} $$ is considered. In this case, the operation is matrix multiplication, and the matrix class ${\mathcal G}$ we may consider as infinite (polyhedron) and finite (its extreme points).
\item[\rm 6.] An $n \times n$ real matrix $\mathbf A$ is called {\it $D$-hyperbolic} if the eigenvalues of ${\mathbf D}{\mathbf A}$ have nonzero real parts for every real nonsingular $n \times n$ diagonal matrix $\mathbf D$. This definition was provided in \cite{AB1}, see also \cite{AB2}, in connection with studying Hopf bifurcation phenomena of linearized system of differential equations. In this case, ${\mathfrak D}$ is the complex plane without imaginary axes, ${\mathcal G}$ is the class of nonsingular diagonal matrices, and $\circ$ is matrix multiplication.
\item[\rm 7.] An $n \times n$ real matrix $\mathbf A$ is called {\it additive $D$-stable} if ${\mathbf D}+{\mathbf A}$ is stable for every positive diagonal matrix $\mathbf D$. According to this, ${\mathfrak D}$ is the open right-hand side of the complex plane, ${\mathcal G}$ is the class of positive diagonal matrices, and $\circ$ is matrix addition. This class was first defined in \cite{CROSS} (referring to the study of diffusion models of biological systems \cite{HAD}) under the name of {\it strong stability}.
\item[\rm 8.] Let, as usual, $[n]$ denotes the set of indices $\{1, \ \ldots, \ n\}$. Given a positive integer $p$, $1 \leq p \leq n$, the set $\alpha = (\alpha_1, \ \ldots \ \alpha_p)$, where each $\alpha_i$ is a nonempty subset of $[n]$, $\alpha_i\bigcap\alpha_j =\emptyset$ and $[n] = \bigcup_i \alpha_i$, is called a {\it partition} of $[n]$. Without loss of generality, we may assume that each $\alpha_i$, $i = 1, \ \ldots, \ p$, consists of contagious indices. A block diagonal matrix $\mathbf H$ with diagonal blocks indexed by $\alpha_1, \ \ldots \ \alpha_p$ is called an {\it $\alpha$-diagonal matrix} (see \cite{HERM}):
    $${\mathbf H} = \diag\{H[\alpha_1], \ \ldots, \ H[\alpha_p]\},$$ where each $H[\alpha_i]$ is a principal submatrix of ${\mathbf H}$ formed by rows and columns with indices from $\alpha_i$, $i = 1, \ \ldots, \ p$.
The following concept was provided in \cite{HERM}: given a partition $\alpha = (\alpha_1, \ \ldots \ \alpha_p)$, an $n \times n$ real matrix $\mathbf A$ is called (multiplicative) {\it $H(\alpha)$-stable} if ${\mathbf H}{\mathbf A}$ is stable for every symmetric positive definite $\alpha$-diagonal matrix $\mathbf H$. Here, the stability region ${\mathfrak D}$ is the open right hand side of the complex plane, the operation $\circ$ is the matrix multiplication, ${\mathcal G}$ is the class of symmetric positive definite $\alpha$-diagonal matrices denoted by $H(\alpha)$. In \cite{GUH1}, in addition, the concept of additive $H(\alpha)$-stability was introduced: given a partition $\alpha = (\alpha_1, \ \ldots \ \alpha_p)$, an $n \times n$ real matrix $\mathbf A$ is called {\it additive $H(\alpha)$-stable} if ${\mathbf H}+{\mathbf A}$ is stable for every symmetric positive definite $\alpha$-diagonal matrix $\mathbf H$. Here, the operation $\circ$ is changed from the matrix multiplication to matrix addition, the stability region ${\mathfrak D}$ and the matrix class ${\mathcal G}$ are as above.
\item[\rm 9.] Given a positive integer $p$, $1 \leq p \leq n$, let $\alpha = (\alpha_1, \ \ldots \ \alpha_p)$ be a partition of $[n]$. A diagonal matrix $\mathbf D$ is called an {\it $\alpha$-scalar matrix} if ${\mathbf D}[\alpha_k]$ is a scalar matrix for every $k = 1, \ \ldots, \ p$, i.e.
    $${\mathbf D} = \diag\{d_{11} {\mathbf I}[\alpha_1], \ \ldots, \ d_{pp} {\mathbf I}[\alpha_p]\}.$$
(Here, as before, ${\mathbf D}[\alpha_k]$ denotes a principal submatrix spanned by rows and columns with indices from $\alpha_k$).
    ${\mathbf D}$ is called a {\it positive $\alpha$-scalar matrix} if, in addition, $d_{ii} > 0$, $i = 1, \ \ldots, \ p$.
    Khalil and Kokotovic introduce the following definition based on the given above matrix class (see \cite{KHAK2}, \cite{KHAK1} ). An $n \times n$ matrix $\mathbf A$ is called {\it $D(\alpha)$-stable} (relative to the partition $\alpha = (\alpha_1, \ \ldots \ \alpha_p)$) if ${\mathbf D}{\mathbf A}$ is stable for every positive $\alpha$-scalar matrix ${\mathbf D}$. (Originally, this property was called {\it block $D$-stability}). The stability region ${\mathfrak D}$ here is the open right-hand side of the complex plane, ${\mathcal G}$ is the class of positive $\alpha$-scalar matrices, $\circ$ is the matrix multiplication (but it is also natural to use matrix addition to introduce the class of {\it additive $D(\alpha)$-stable} matrices).
 \item[\rm 10.]   The following definition was provided in \cite{KU1}: given a positive diagonal matrix ${\mathbf D} = \diag\{d_{11}, \ \ldots, \ d_{nn}\}$ and a permutation $\theta = (\theta(1), \ \ldots, \ \theta(n))$ of the set of indices $[n]:= \{1, \ \ldots, \ n\}$, we call the matrix ${\mathbf D}$ {\it ordered with respect to $\theta$}, or {\it $\theta$-ordered}, if it satisfies the inequalities $$d_{\theta(i)\theta(i)} \geq d_{\theta(i+1)\theta(i+1)}, \qquad i = 1, \ \ldots, \ n-1. $$ A matrix $\mathbf A$ is called {\it $D$-stable with respect to the order $\theta$}, or {\it $D_\theta$-stable}, if the matrix ${\mathbf D}{\mathbf A}$ is positive stable for every $\theta$-ordered positive diagonal matrix $\mathbf D$.
     Moreover, if we consider a union of $k$ classes, defined by the permutations $\theta_1, \ \ldots, \ \theta_k$, respectively, we denote it ${\mathcal G}_{\theta_1, \ldots, \theta_k}$.
 \item[\rm 11.] A matrix $\mathbf A$ is called {\it $D$-stable with respect to the set $\Theta \subset {\mathcal M}^{n \times n}$} if ${\mathbf D}{\mathbf A}$ is stable for every matrix ${\mathbf D} \in \Theta$. This approach is used by Kosov for studying multiplicative $D$-stability (see \cite{KOS}). Note, that this is not a proper decomposition of the set of all positive diagonal matrices.
For studying additive $D$-stability, the following set is considered (see \cite{KOS}, \cite{RS1}):
$$\Theta_0 = \prod(0, \ d_{ii}^{max})= $$ $$ \diag (d_{ii}, \ 0 < d_{ii} < d_{ii}^{max} < +\infty, \ \ i = 1, \ 2, \ \ldots \ n). $$
\item[\rm 12.] A real matrix ${\mathbf A}$ is called {\it Hadamard $H$-stable} if ${\mathbf H}\circ{\mathbf A}$ is stable for every symmetric positive definite matrix ${\mathbf H} \in {\mathcal M}^{n \times n}$. Here the domain ${\mathfrak D}$ is the open right-hand side of the complex plane, ${\mathcal G}$ is the class of symmetric positive definite matrices, $\circ$ is the Hadamard (entry-wise) matrix multiplication. This definition was introduced by Johnson (see \cite{JOHN3}, p. 304) and was called by the author {\it Schur stability}. But, since the term "Schur stable" is already reserved for matrices which spectral radius is less than 1, we refer to this property as to {\it Hadamard $H$-stability}.
\item[\rm 13.] The following matrix classes were introduced in \cite{JOHD} in order to "interpolate" the matrix properties and the properties of sign pattern classes, in particular, the classes of $D$-stable and sign-stable matrices. An $n \times n$ real matrix $\mathbf A$ is called {\it $B_k$-stable} (to belong to the class $B_k$) if the Hadamard product ${\mathbf B}\circ{\mathbf A}$ is positive stable for every entry-wise positive matrix ${\mathbf B} \in {\mathcal M}^{n \times n}$ with ${\rm rank}({\mathbf B}) \leq k$. Here, the domain ${\mathfrak D}$ is the open right-hand side of the complex plane, ${\mathcal G}$ is the class of entry-wise positive matrices of finite rank $k$ and $\circ$ is the Hadamard matrix multiplication. By varying the class ${\mathcal G}$, the authors also introduce the class of $B_k^+$-stable matrices: an $n \times n$ real matrix $\mathbf A$ is called {\it $B_k^+$-stable} (to belong to the class $B_k^+$) if the Hadamard product ${\mathbf B}\circ{\mathbf A}$ is positive stable for every matrix ${\mathbf B} \in {\mathcal M}^{n \times n}$ such that $${\mathbf B} = {\mathbf B}_1 + {\mathbf B}_2 + \ldots + {\mathbf B}_k,$$ where each ${\mathbf B}_i$ is entry-wise positive, ${\rm rank}({\mathbf B}_i) = 1$, $i = 1, \ 2, \ \ldots, \ k$.
 \item[\rm 14.] Changing the stability region ${\mathfrak D}$ to ${\mathbb C} \setminus \{0\}$, the same kind of interpolation between nonsingular and sign nonsingular matrices was introduced in \cite{JOHD} (see \cite{JOHD}, p. 368). An $n \times n$ matrix $\mathbf A$ is called {\it $B_k$-nonsingular} if the Hadamard product ${\mathbf B}\circ{\mathbf A}$ is nonsingular for every entry-wise positive matrix ${\mathbf B} \in {\mathcal M}^{n \times n}$ with ${\rm rank}({\mathbf B}) \leq k$ (in \cite{JOHD}, this matrix class is denoted by $L_{n,k}$). For strong forms of nonsingularity, see also \cite{DJD}.
\item[\rm 15.]  Given an $n \times n$ matrix $\mathbf A$, consider its finite-rank perturbation of the form
\begin{equation}\label{pert} \widetilde{{\mathbf A}} = {\mathbf A} + {\mathbf B}, \end{equation}
where ${\mathbf B} \in {\mathcal M}^{n \times n}$ with ${\rm rank}({\mathbf B}) \leq k$, $k = 1, \ \ldots, \ n$. If ${\rm rank}({\mathbf B}) = 1$, Equality \eqref{pert} may be written as
$$ \widetilde{{\mathbf A}} = {\mathbf A} + x \otimes y,$$
where $x, \ y \in {\mathbb R}^n$. The following problem studied by Barkovsky (see \cite{BARK}, \cite{BAY}): given an $n \times n$ matrix $\mathbf A$ and two vectors $x, \ y \in {\mathbb R}^n$, when all the matrices $\widetilde{{\mathbf A}}_{\tau} = {\mathbf A} +\tau (x \otimes y)$ have real spectra? This problem can be considered as establishing $({\mathfrak D},{\mathcal G},\circ)$-stability, where the stability region $\mathfrak D$ is the real axes, ${\mathcal G}$ is a parametric rank-one matrix family of the form $\{\tau (x\otimes y)\}_{\tau \in {\mathbb R}}$, the operation $\circ$ is matrix addition. For the problems of this type, see also \cite{FGR}, \cite{BIR}.
\end{enumerate}

The following generalization of $D$-stability for two fixed stability regions: the open right-hand side of the complex plane and the interior of the unit disk was defined in \cite{CADHJ} (see \cite{CADHJ}, p. 152). By fixing the operation $\circ$ as matrix multiplication and varying the class ${\mathcal G} \subset {\mathcal M}^{n \times n}$, the following definitions were obtained: an $n \times n$ real matrix $\mathbf A$ is called {\it ${\mathcal G}$-stable} ({\it ${\mathcal G}$-convergent}) if ${\mathbf G}{\mathbf A}$ is positive stable (respectively, convergent) for every matrix ${\mathbf G}$ from the selected matrix class ${\mathcal G}$. The concept of "set product", when a matrix class ${\mathcal G}$ is multiplied by a given matrix $\mathbf A$ is further studied in \cite{BHK3} and \cite{CALN}.

\subsection{Stability regions}
The theory of stability in some generalized regions known as {\it root clustering} was mainly developed in 1980th by Gutman and Jury (see \cite{GUJU}) and continued in \cite{GUT}, \cite{GUT2}, \cite{JU2}, where a lot of special classes of stability regions were analyzed. Here, we refer to \cite{HHP}, \cite{LLK}, \cite{MEME}, \cite{MAO} (disk regions), \cite{ANB}, \cite{BICJ1}, \cite{ANB}, \cite{BICJ2}, \cite{BICJ3}, \cite{MOMOK}, \cite{SGH} (sector regions). The methods of studying $({\mathfrak D},{\mathcal G},\circ)$-stability are mainly defined by the methods of stability region description and of studying $\mathfrak D$-stable matrices.

The approach we will use later is based on the necessary and sufficient condition of matrix stability, proved by Lyapunov (see, for example, \cite{BELL}, \cite{GANT2}, for exact formulation see \cite{GANT}, \cite{HER1}, p. 164, Theorem 2.4).
    \begin{theorem}[Lyapynov] An $n \times n$ matrix $\mathbf A$ is (positive) stable if and only if there exists a symmetric positive definite matrix $\mathbf H$ such that the matrix
    $${\mathbf H}{\mathbf A} + {\mathbf A}^{T}{\mathbf H}$$
    is positive definite.
    \end{theorem}
Equivalently, we analyze the solvability of the Lyapunov equation
\begin{equation}\label{lyap}
{\mathbf H}{\mathbf A} + {\mathbf A}^{T}{\mathbf H} = {\mathbf W},
\end{equation}
where ${\mathbf W}$ is a symmetric positive definite matrix, in the class $\mathcal{ H}$ of symmetric positive definite matrices.

For a class of more general stability regions $\mathfrak D$, we need the following generalization of Lyapunov theorem, obtained by Hill (see \cite{HIL}, also \cite{WIM}, p. 140, Theorem 1).

\begin{theorem}\label{TGenLyap} If an $n \times n$ matrix $\mathbf A$ satisfies the matrix equation
\begin{equation}\label{GenLyap} \sum_{i,j = 0}^{n-1}c_{ij}({\mathbf A}^T)^i{\mathbf H}{\mathbf A}^j = {\mathbf W}, \qquad c_{ij} = c_{ji}
\end{equation}
with a symmetric positive definite matrix $\mathbf H$, then
\begin{enumerate}
\item[\rm 1.] ${\mathbf W}$ is a symmetric positive definite matrix implies $$f(\lambda) := \sum_{i,j = 0}^{n-1}c_{ij}\overline{\lambda}^i\lambda^j > 0;$$
\item[\rm 2.] ${\mathbf W}$ is a symmetric positive semidefinite matrix implies $$f(\lambda) := \sum_{i,j = 0}^{n-1}c_{ij}\overline{\lambda}^i\lambda^j \geq 0;$$
\item[\rm 3.] ${\mathbf W} = {\mathbf 0}$ implies $$f(\lambda) := \sum_{i,j = 0}^{n-1}c_{ij}\overline{\lambda}^i\lambda^j = 0.$$
\end{enumerate}
\end{theorem}

Consider the special cases of Theorem \ref{TGenLyap} describing the most used stability regions:

\begin{enumerate}
\item[\rm 1.] {\bf Case I. ${\mathfrak D} = \{\lambda \in {\mathbb C}: {\rm Re}(\lambda)> 0\}$.} Here, we deal with the classical Lyapunov equation \eqref{lyap}.
\item[\rm 2.] {\bf Case II. ${\mathfrak D} = \{\lambda \in {\mathbb C}: {\rm Re}(\lambda) \neq 0\}$}.
    Here, instead of Lyapynov theorem, we deal with the following theorem proved in \cite[p. 76]{OSS} (see \cite{OSS}, p. 76, Theorem 1).
    \begin{theorem} An $n \times n$ matrix $\mathbf A$ has no pure imaginary eigenvalues (i.e. with zero real parts) if and only if there exists a symmetric matrix $\mathbf H$ such that the matrix
    $${\mathbf W}:={\mathbf H}{\mathbf A} + {\mathbf A}^{T}{\mathbf H}$$
    is positive definite. Then we have ${\rm In}({\mathbf H}) = {\rm In}({\mathbf A})$.
    \end{theorem}

\item[\rm 3.] {\bf Case III. ${\mathfrak D} = \{\lambda \in {\mathbb C}: |\lambda| < 1\}$}. Here, an analougous statement will be used (see \cite{STE}, also \cite{WIM}).
    \begin{theorem}[Stein] An $n \times n$ matrix $\mathbf A$ is Schur stable if and only if there exists a symmetric positive definite matrix $\mathbf H$ such that the matrix
    \begin{equation}\label{ste}{\mathbf W}:={\mathbf H} - {\mathbf A}^{T}{\mathbf H}{\mathbf A}\end{equation}
    is positive definite.
    \end{theorem}
\end{enumerate}

\subsection{Matrix classes and their properties} Here, we collect and analyze the most studied cases of matrix classes $\mathcal G$. We are especially interested in the following basic facts.
\begin{enumerate}
\item[\rm 1.] The inclusion relations between the studied matrix classes.
\item[\rm 2.] For a given group operation $\circ$ on ${\mathcal M}^{n \times n}$, if ${\mathcal G}$ is closed with respect to this operation, moreover, if $({\mathcal G}, \circ)$ form a subgroup.
\item[\rm 3.] Commutators of the class ${\mathcal G}$ and transformations that leave this class invariant.
\end{enumerate}

 Now let us consider the following matrix classes.

 \begin{enumerate}
\item[\rm 1.] Class ${\mathcal S}$ of symmetric matrices from ${\mathcal M}^{n \times n}$. This matrix class, as well as all its subclasses, is closed with respect to matrix transposition. Class ${\mathcal S}$ equipped with the operation of matrix addition forms a group, nonsingular symmetric matrices form a group with respect to matrix multiplication, and matrices without zero entries form a group with respect to Hadamard matrix multiplication.
To study commutators, we need the following lemma (see \cite{HOJ}, p. 172).
\begin{lemma} Let $\mathbf A$, ${\mathbf B}$ be symmetric matrices. Then ${\mathbf A}{\mathbf B}$ is also symmetric if and only if $\mathbf A$ and $\mathbf B$ commute.
\end{lemma}
\item[\rm 2.] Class ${\mathcal H}$ of symmetric positive definite matrices. Here, we mention the following characterization (see \cite{BHAT2}, p. 2): {\it a matrix ${\mathbf A} \in {\mathcal M}^{n \times n}$ is symmetric positive definite if and only if ${\mathbf A} = {\mathbf B}^T{\mathbf B}$ for some matrix ${\mathbf B} \in {\mathcal M}^{n \times n}$}. The class of symmetric positive definite matrices is closed under Hadamard multiplication (first proved in \cite{SHU}, see also \cite{BHAT2}) and under matrix addition. However, for ${\mathbf A}, {\mathbf B} \in {\mathcal H}$, the usual matrix product ${\mathbf A}{\mathbf B}$ belongs to ${\mathcal H}$ if and only if $\mathbf A$ and $\mathbf B$ commute (\cite{BHAT2}). Class ${\mathcal H}$ is also closed with respect to multiplicative inverse.
\item[\rm 3.] Class ${\mathcal H}_{\alpha}$ of symmetric $\alpha$-diagonal matrices, for a given partition $\alpha = (\alpha_1, \ \ldots \ \alpha_p)$ of $[n]$. Recall, that an $n \times n$ matrix $\mathbf H$ is called an $\alpha$-diagonal matrix if its principal submatrices ${\mathbf H}[\alpha_i]$  (formed by rows and columns with indices from $\alpha_i$, $i = 1, \ \ldots, \ p$) are nonzero, while the rest of the entries is zero. This class forms a group with respect to matrix addition.
\item[\rm 4.] Class $\mathcal D$ of diagonal matrices. This class form a group under operation of matrix addition. Nonsingular diagonal matrices also form an abelian group with respect to matrix multiplication.
\item[\rm 5.] Sign pattern classes ${\mathcal D}_S$. First, define a {\it sign pattern} ${\rm Sign}({\mathbf D})$ of a diagonal matrix $\mathbf D$ as follows: $${\rm Sign}({\mathbf D}) := {\rm diag}\{{\rm sign}(d_{11}), \ \ldots, \ {\rm sign}(d_{nn})\}.$$
Two diagonal matrices ${\mathbf D}_1$ and ${\mathbf D}_2$ are said to belong to the same sign pattern class if ${\rm Sign}({\mathbf D}_1) = {\rm Sign}({\mathbf D}_2)$. For a given sign pattern $S$, define ${\mathcal D}_S$ as a sign pattern class of diagonal matrices.
 The set of all sign pattern classes covers the set of all diagonal matrices. So we have a proper decomposition
 $${\mathcal G} = \bigcup_{S}{\mathcal G}(S).$$
 Sign pattern classes are studied in connection with matrix inertia properties, $D$-hyperbolicity and Schur $D$-stability.
 \item[\rm 6.] Class ${\mathcal D}^+$ of positive diagonal matrices. This class is closed with respect to matrix addition forms abelian group matrix multiplications.
 \item[\rm 7.] Class ${\mathcal D}_{\alpha}$ of $\alpha$-scalar matrices (resp. ${\mathcal D}^+_{\alpha}$ of positive $\alpha$-scalar matrices). Recall that, for a given partition $\alpha=(\alpha_1, \ \ldots \ \alpha_p)$ of $[n]$, $1 \leq p \leq n$, a diagonal matrix $\mathbf D$ is called an {\it $\alpha$-scalar matrix} if ${\mathbf D}[\alpha_k]$ is a scalar matrix for every $k = 1, \ \ldots, \ p$, i.e.
    $${\mathbf D} = \diag\{d_{11} {\mathbf I}_{\alpha_1}, \ \ldots, \ d_{pp} {\mathbf I}_{\alpha_p}\}.$$
    ${\mathbf D}$ is called a {\it positive $\alpha$-scalar matrix} if, in addition, $d_{ii} > 0$, $i = 1, \ \ldots, \ p$.
    For studying this matrix class, see \cite{HERM}, \cite{WAN}, \cite{KHAK2} and \cite{AB1}, \cite{AB2}.
    For a fixed $\alpha$, the class ${\mathcal D}_{\alpha}$ is closed with respect to matrix addition and forms abelian group with respect to matrix multiplication.
\item[\rm 8.] Class ${\mathcal D}_{\theta}$ of positive diagonal matrices ordered with respect to a given permutation $\theta \in \Theta_{[n]}$ and a union ${\mathcal D}_{\theta_1, \ldots, \theta_k}$ of $k$ classes, defined by the permutations $\theta_1, \ \ldots, \ \theta_k$ (Recall that a positive diagonal matrix ${\mathbf D} = \diag\{d_{11}, \ \ldots, \ d_{nn}\}$ is called {\it ordered with respect to a permutation $\theta = (\theta(1), \ \ldots, \ \theta(n))$ of $[n]$}, or {\it $\theta$-ordered}, if it satisfies the inequalities $$d_{\theta(i)\theta(i)} \geq d_{\theta(i+1)\theta(i+1)}, \qquad i = 1, \ \ldots, \ n-1. $$
     In what follows, $\Theta_{[n]}$ denotes, as usual, the set of all the permutations of $[n]$. Obviously, $${\mathcal D} = \bigcup_{\theta \in \Theta_{[n]}}{\mathcal D}_{\theta}.$$
     This class is closed with respect to matrix addition and matrix multiplication, however, does not contain multiple inverses.
\item[\rm 9.] Class ${\mathcal D}_{\Theta}$ of diagonal matrices satisfying the inequalities
$$\Theta = \prod(d_{ii}^{min}, \ d_{ii}^{max})= $$ $$ \diag (d_{ii}, \ 0 < d_{ii}^{min} < d_{ii} < d_{ii}^{max} < +\infty, \ \ i = 1, \ 2, \ \ldots \ n); $$
and class ${\mathcal D}_{\Theta_0}$, where $$\Theta_0 = \prod(0, \ d_{ii}^{max})= $$ $$ \diag (d_{ii}, \ 0 < d_{ii} < d_{ii}^{max} < +\infty, \ \ i = 1, \ 2, \ \ldots \ n). $$
\item[\rm 10.] Class ${\mathcal D}_V$ of vertex diagonal matrices. Recall that a real diagonal matrix $\mathbf D$ is called {\it vertex diagonal} if $|{\mathbf D}|=1 $, i.e. $|d_{ii}| = 1$ for any $i = 1, \ \ldots, \ n$ (\cite{BHK}).
\item[\rm 11.]
 Now let us consider the following characterizations of a $\theta$-ordered matrix $D$:
 $$d_{max}({\mathbf D}) := \max_i\dfrac{d_{\theta(i)\theta(i)}}{d_{\theta(i+1)\theta(i+1)}}, \qquad i = 1, \ \ldots, \ n.$$
$$d_{min}({\mathbf D}) := \min_i\dfrac{d_{\theta(i)\theta(i)}}{d_{\theta(i+1)\theta(i+1)}}, \qquad i = 1, \ \ldots, \ n.$$
According to the definition, $1 \leq d_{min}({\mathbf D}) \leq d_{max}({\mathbf D})$ for every positive diagonal matrix $\mathbf D$.
\end{enumerate}
Here, we have the following chains of inclusions:
\begin{equation}\label{chain}
{\mathcal D}_{\alpha}^+ \subset {\mathcal D}^+ \subset {\mathcal H}_{\alpha} \subset {\mathcal H}.
\end{equation}

Now let us consider the following pairwise commuting subclasses of ${\mathcal H}$:
\begin{enumerate}
\item[\rm 1.] $({\mathcal H}, {\mathcal I})$, where the class ${\mathcal I}$ consists of the only one identity matrix $\mathbf I$.
\item[\rm 2.] $({\mathcal H}_{\alpha}, {\mathcal D}_{\alpha})$ ($\alpha$-scalar diagonal matrices commute with $\alpha$-block symmetric matrices).
\item[\rm 3.] $({\mathcal D}, {\mathcal D})$ (diagonal matrices commute within themselves).
\end{enumerate}

\subsection{Binary operations}
\paragraph{Basic definitions and properties}
Let us recall the following definitions and properties we will use later.
Consider a binary operation $\circ$ on a matrix class ${\mathcal G}_0 \subseteq {\mathcal M}^{n \times n}$:
$$\circ:{\mathcal G}_0 \times {\mathcal G}_0 \rightarrow {\mathcal M}^{n \times n}.$$For the further study, it would be convenient to assume that the class ${\mathcal G}$ belongs to ${\mathcal G}_0$ to avoid the question of spreading the binary operation $\circ$ to the matrices from ${\mathcal G}$. Let us mention the following operation properties (see, for example, \cite{CURT}).
\begin{enumerate}
\item[\rm 1.] Associativity
$${\mathbf A}\circ({\mathbf B}\circ{\mathbf C}) = ({\mathbf A}\circ{\mathbf B})\circ{\mathbf C}$$
for every ${\mathbf A}, {\mathbf B}, {\mathbf C} \in {\mathcal G}_0$.
\item[\rm 2.] There exists an {\it identity element} ${\mathbf L} \in {\mathcal G}_0$:
$${\mathbf L}\circ{\mathbf A} = {\mathbf A}\circ{\mathbf L} = {\mathbf A}$$
for every ${\mathbf A} \in {\mathcal G}_0$.
\item[\rm 3.] There exist {\it inverses}: for every ${\mathbf A} \in {\mathcal G}_0$, there is ${\mathbf A}^{-1} \in {\mathcal G}_0$ such that
$${\mathbf A}\circ{\mathbf A}^{-1} = {\mathbf A}^{-1}\circ{\mathbf A} = {\mathbf L}.$$
\item[\rm 4.] Commutativity
$${\mathbf A}\circ{\mathbf B} = {\mathbf B}\circ{\mathbf A}$$
for every ${\mathbf A}, {\mathbf B} \in {\mathcal G}_0$.
\end{enumerate}

A {\it group} $({\mathcal G}_0, \circ)$ is a set ${\mathcal G}_0$ equipped with a binary operation $\circ$, which satisfies Properties 1-3. If, in addition, the operation $\circ$ satisfies Property 4, a group $({\mathcal G}_0, \circ)$ is called {\it abelian}. If ${\mathcal G} \subset {\mathcal G_0} \subseteq {\mathcal M}^{n \times n}$ and is closed with respect to $\circ$ then $({\mathcal G}, \circ)$ is a subgroup of $({\mathcal G}_0, \circ)$ if $\circ$ satisfies Properties 1-3 on $\mathcal G$. Any matrix property, that is preserved under operation $\circ$, the identity element and inverses with respect to $\circ$ can form a subgroup with respect to $\circ$.

 A group $({\mathcal G}_0, \circ)$ is called {\it topological} if it is a topological space and the group operation $\circ$ is continuous in this topological space. For the definition and theory of topological groups, we refer to \cite{PONT}, for more detailed study of the question see \cite{ART}). In some cases, we will also assume that the class ${\mathcal G}$ forms a subgroup of the topological group ${\mathcal G}_0$ (i.e. a subgroup, which is a closed subspace in the topological space ${\mathcal G}_0$).

We may also consider more theoretical examples, which arises mostly in control theory, i.e. Lyapunov operator (see \cite{BHAT2}) and its generalizations, block Hadamard product (see \cite{HMN}, \cite{CHO}), Redheffer product (see \cite{TIM}), Hurwitz product (for the definition see, for example, \cite{ALAL}), the max-algebra operations (see \cite{BUT}), sub-direct sums (see \cite{FAJ}).

\paragraph{Relations to matrix addition and matrix multiplication} 

First, let us consider operation of {\bf matrix addition}. Later, we need the following cases.
\begin{enumerate}
\item[\rm 1.] The operation $\circ$ and matrix addition $+$ are connected with the rule of associativity
\begin{equation}\label{+ass}({\mathbf A} \circ {\mathbf B}) + {\mathbf C} = {\mathbf A} \circ ({\mathbf B} + {\mathbf C}).\end{equation}
For example, $\circ$ is also matrix addition.
\item[\rm 2.] The operation $\circ$ is distributive over $+$
\begin{equation}\label{+dist1}({\mathbf A} + {\mathbf B})\circ{\mathbf C} = ({\mathbf A} \circ {\mathbf C}) + ({\mathbf B}\circ{\mathbf C}). \end{equation}
Here, we consider the operations of usual and Hadamard matrix multiplication.
\item[\rm 3.] The operation $+$ is distributive over $\circ$
\begin{equation}\label{+dist2}{\mathbf A} +({\mathbf B}\circ {\mathbf C}) = ({\mathbf A} +{\mathbf B}) \circ ({\mathbf A} +{\mathbf C}).\end{equation}
As an example, we consider the operation of entry-wise maximum.
\end{enumerate}
Note, that for an arbitrary operation $\circ$ none of the above formulae may hold.

{\bf Scalar multiplication}
\begin{enumerate}
\item[\rm 1.] The operation of multiplication by a scalar $\alpha \in {\mathbb R}$ is connected to the operation $\circ$ by the rules of associativity and commutativity:
\begin{equation}\label{scaas} \alpha({\mathbf A}\circ{\mathbf B}) = (\alpha{\mathbf A})\circ{\mathbf B} = {\mathbf A}\circ(\alpha{\mathbf B})\end{equation}
for every ${\mathbf A}, {\mathbf B} \in {\mathcal M}^{n \times n}$.

{\bf Examples.} Matrix multiplication, Hadamard matrix multiplication.
\item[\rm 2.] The operation of scalar multiplication is connected to the operation $\circ$ by the rule of distributivity:
\begin{equation}\label{scadist}\alpha({\mathbf A}\circ{\mathbf B}) = (\alpha{\mathbf A})\circ(\alpha{\mathbf B})\end{equation}
for every ${\mathbf A}, {\mathbf B} \in {\mathcal M}^{n \times n}$.

{\bf Examples.} Matrix addition, matrix maximum.
\end{enumerate}

{\bf Matrix multiplication.}
Now let us consider the relations between $\circ$ and matrix multiplication. Later, we consider one of the following cases:
\begin{enumerate}
\item[\rm 1.] Associativity:
\begin{equation}\label{*ass}{\mathbf A}\circ({\mathbf B}{\mathbf C}) = ({\mathbf A}{\mathbf B})\circ{\mathbf C}\end{equation}
Examples: matrix multiplication.
\item[\rm 2.] Distributivity (the operation $\circ$ is distributive over matrix multiplication):
\begin{equation}\label{*dist1}{\mathbf A}\circ({\mathbf B}{\mathbf C}) = ({\mathbf A}\circ{\mathbf B})({\mathbf A}\circ{\mathbf C})\end{equation}
\item[\rm 3.] Distributivity (matrix multiplication is distributive over $\circ$):
\begin{equation}\label{*dist2}{\mathbf A}({\mathbf B}\circ{\mathbf C}) = ({\mathbf A}{\mathbf B})\circ({\mathbf A}{\mathbf C})\end{equation}
Examples: matrix addition.
\end{enumerate}

\section{Elementary properties of $({\mathfrak D}, \ {\mathcal G}, \ \circ)$-stability classes}
\subsection{Inclusion relations and topological properties} Here, we collect basic statements describing the class of $({\mathfrak D}, \ {\mathcal G}, \ \circ)$-stable matrices.

First, let us mention a topological property, which shows if the definition of $({\mathfrak D}, \ {\mathcal G}, \ \circ)$-stability is meaningful, i.e. if the defined class of $({\mathfrak D}, \ {\mathcal G}, \ \circ)$-stable matrices is nonempty.

\begin{theorem} Given a bounded (in absolute value) stability region ${\mathfrak D} \subseteq {\mathbb C}$, a matrix class ${\mathcal G} \subset {\mathcal M}^{n \times n}$ and a continuous binary operation $\circ$ on ${\mathcal M}^{n \times n} \times {\mathcal M}^{n \times n}$. Then, for the class of $({\mathfrak D},{\mathcal G},\circ)$-stable matrices to be nonempty it is necessary the class ${\mathcal G}$ be bounded in ${\mathcal M}^{n \times n}$ with respect to any matrix norm.
\end{theorem}
{\bf Proof.} Let the matrix class $\mathcal G$ be unbounded, i.e. there exists a sequence $\{{\mathbf G}_i\}_{i=1}^{\infty}$ from ${\mathcal G}$ with $\|{\mathbf G}_i\| \rightarrow \infty$. Let there exist at least one $({\mathfrak D}, \ {\mathcal G}, \ \circ)$-stable matrix $\mathbf A$. By definition,  $\sigma({\mathbf G}\circ {\mathbf A}) \subset {\mathfrak D}$ for all ${\mathbf G} \in {\mathcal G}$. Since ${\mathfrak D}$ is bounded in $\mathbb C$, there is a positive value $R$ such that $|\lambda| \leq R$ for all $\lambda \in {\mathfrak D}$. Thus the spectral radius $\rho({\mathbf G}\circ {\mathbf A}) \leq R$ for all ${\mathbf G} \in {\mathcal G}$. By the continuity of the operation $\circ$, we have
 $\|{\mathbf G}_i\circ {\mathbf A}\|\rightarrow\infty $ as $\|{\mathbf G}_i\| \rightarrow \infty$ for any fixed $\mathbf A$, which implies
 $\rho({\mathbf G}_i\circ {\mathbf A}) \rightarrow \infty$ as well. Thus we came to a contradiction. $\square$

{\bf Example.} Consider any bounded stability region ${\mathfrak D} \subseteq {\mathbb C}$ and ${\mathcal M}^{n \times n}$ equipped with the operations of matrix addition and matrix multiplication. Then the inequalities
$$\|{\mathbf A} + {\mathbf G}\| \geq \|{\mathbf G} \| - \|{\mathbf A}\|;$$
$$\|{\mathbf A}{\mathbf G}\| \geq \|{\mathbf G}\|\frac{1}{\|{\mathbf A}^{-1}\|}$$
imply the boundness of the matrix class $\mathcal G$.

The next results deal with basic inclusion relations between different classes of $({\mathfrak D},{\mathcal G},\circ)$-stable matrices.

\begin{theorem}\label{mainob1} Let ${\mathfrak D}\subset {\mathbb C}$ be a stability region of the complex plane, ${\mathcal G} \subset {\mathcal M}^{n \times n}$ be a matrix class and $\circ$ be a binary operation on ${\mathcal M}^{n \times n}$. Then
     \begin{enumerate}
\item[\rm 1.] for any subset ${\mathfrak D}_1$ of the complex plane such that ${\mathfrak D}_1 \subseteq {\mathfrak D}$, the class of $({\mathfrak D}_1, \ {\mathcal G}, \ \circ)$-stable matrices belongs to the class of $({\mathfrak D}, \ {\mathcal G}, \ \circ)$-stable matrices.
\item[\rm 2.] for any matrix class ${\mathcal G}_1$ such that ${\mathcal G}_1 \subseteq {\mathcal G}$, conversely, the class of $({\mathfrak D}, \ {\mathcal G}, \ \circ)$-stable matrices belongs to the class of $({\mathfrak D}, \ {\mathcal G}_1, \ \circ)$-stable matrices. In particular, if $${\mathcal G} = \bigcup_{i \in I}{\mathcal G}_{i},$$
a matrix $\mathbf A$ is $({\mathfrak D}, \ {\mathcal G}, \ \circ)$-stable if and only if it is $({\mathfrak D}, \ {\mathcal G}_i, \ \circ)$-stable for any $i \in I$.
\end{enumerate}
     \end{theorem}
  {\bf Proof.} The proof obviously follows from the definition of $({\mathfrak D}, \ {\mathcal G}_1, \ \circ)$-stability.

{\bf Example 1.} Let ${\mathcal G}$ be the class of positive diagonal matrices, $\circ$ be matrix multiplication. The following inclusion relations between $({\mathfrak D}, \ {\mathcal G}, \ \circ)$-stability classes are based on the inclusion relations between the corresponding stability regions (the positive direction of the real axes belongs to the open right half plane of the complex plane which belongs to the complex plane without the imaginary axes):
$$\mbox{$D$-positive matrices} \subset \mbox{$D$-stable matrices} \subset \mbox{$D$-hyperbolic matrices}.$$

{\bf Example 2.} The following relations are based on the inclusion chain \ref{chain} between matrix classes.
$$\mbox{$H$-stable matrices} \subset \mbox{$H_{\alpha}$-stable matrices} \subset \mbox{$D$-stable matrices} \subset $$ $$\mbox{$D_{\alpha}$-stable matrices} \subset \mbox{stable matrices}.$$
The relation $$\mbox{$H$-stable matrices} \subset \mbox{$D$-stable matrices}$$ is pointed out in \cite{AM}.

{\bf Example 3.} Let $\alpha$ and $\beta$ be two partitions of the set $[n]$, such that $\beta \subseteq \alpha$. Then
$$\mbox{$H_{\beta}$-stable matrices} \subset \mbox{$H_{\alpha}$-stable matrices} $$
and
$$\mbox{$D_{\alpha}$-stable matrices} \subset \mbox{$D_{\beta}$-stable matrices}. $$

The inclusion $$\mbox{additive $H_{\beta}$-stable matrices} \subset \mbox{additive $H_{\alpha}$-stable matrices} $$
was pointed out in \cite{GUH1} (see \cite{GUH1}, p. 327, Theorem 2.1).

{\bf Example 4.}
The class of Schur $D$-stable matrices belongs to the class of vertex $D$-stable matrices. In some cases ($n =3$ (see \cite{FL1}, p. 20, Theorem 2.7), tridiagonal matrices (see \cite {FL2}, p. 46, Theorem 4.2), some others) these classes coincide.

{\bf Example 5.}
$$\mbox{$D_{\theta}$-stable matrices} \subset \mbox{$D$-stable matrices},$$
for any permutation $\theta$ of the set $[n]$. Let us consider the decompositions of the class of $D$-stable matrices. The following statement was proven in \cite{KU1}.
 \begin{lemma}\label{Dst}
A matrix $\mathbf A$ is $D$-stable if and only if it is $D_\theta$-stable for any $\theta \in \Theta_{[n]}$.
\end{lemma}

Now we make the most common observation which may be used describing the class of $({\mathfrak D},{\mathcal G},\circ)$-stable matrices.
 \begin{theorem}\label{mob} Given a stability region ${\mathfrak D} \subset {\mathbb C}$, a matrix class ${\mathcal G}\subset{\mathcal M}^{n \times n}$ and a binary operation $\circ$ on ${\mathcal M}^{n \times n}$,
any condition on matrices which implies ${\mathfrak D}$-stability and which is preserved under ${\mathbf G}\circ(\cdot)$ for any $\mathbf G$ from the class ${\mathcal G}$, is sufficient for $({\mathfrak D},{\mathcal G},\circ)$-stability.
\end{theorem}
The proof is obvious.
 This is a generalization of a simple observation, made by Johnson for the case of multiplicative $D$-stability (see \cite{JOHN1}, p. 54, Observation (i)).

However, in each special case, there may be other classes, not covered by this reasoning.

Finally, let us consider the following result on the belonging to the class of $({\mathfrak D},{\mathcal G},\circ)$-stable matrices and the connection to ${\mathfrak D}$-stability.

\begin{theorem}\label{Dstab} Let $\mathfrak D$ be an arbitrary stability region and $\circ$ be a group operation on ${\mathcal M}^{n \times n}$. Then
 \begin{enumerate}
 \item[\rm 1.] If the matrix class $\mathcal G \subset {\mathcal M}^{n \times n}$ is closed with respect to the operation $\circ$, then $\mathbf A$ is $({\mathfrak D},{\mathcal G},\circ)$-stable implies ${\mathbf G} \circ {\mathbf A}$ is $({\mathfrak D},{\mathcal G},\circ)$-stable for any ${\mathbf A} \in {\mathcal M}^{n \times n}$ and any ${\mathbf G} \in {\mathcal G}$.
 \item[\rm 2.] If the matrix class $\mathcal G$ includes the identity element $\mathbf L$ (with respect to the operation $\circ$), then every $({\mathfrak D},{\mathcal G}, \circ)$-stable matrix is ${\mathfrak D}$-stable.
\item[\rm 3.] If $\mathcal G$ forms a subgroup with respect to the operation $\circ$, then ${\mathbf G}_0\circ{\mathbf A}$ is $({\mathfrak D},{\mathcal G},\circ)$-stable for at least one matrix ${\mathbf G}_0 \in {\mathcal G}$ implies ${\mathbf A}$ is $\mathfrak D$-stable and $({\mathfrak D},{\mathcal G},\circ)$-stable.
\end{enumerate}
\end{theorem}
{\bf Proof}. \begin{enumerate}
 \item[\rm 1.] Let $\mathbf A$ be $({\mathfrak D},{\mathcal G},\circ)$-stable. Consider ${\mathbf G}\circ{\mathbf A}$, for an arbitrary ${\mathbf G} \in {\mathcal G}$. Then, taking an arbitrary ${\mathbf G}_0 \in {\mathcal G}$, we obtain
$${\mathbf G}_0\circ({\mathbf G}\circ{\mathbf A}) = [\mbox{associativity}] = ({\mathbf G}_0\circ{\mathbf G})\circ{\mathbf A} =$$
$$ = {\mathbf G_1}\circ{\mathbf A},$$
where ${\mathbf G_1}:={\mathbf G}_0\circ{\mathbf G} \in {\mathcal G}$ due to its closeness. Thus $\sigma({\mathbf G}_0\circ({\mathbf G}\circ{\mathbf A})) = \sigma({\mathbf G_1}\circ{\mathbf A}) \subset {\mathfrak D}$ for any ${\mathbf G} \in {\mathcal G}$.
\item[\rm 2.] Let ${\mathbf A}$ be a $({\mathfrak D},{\mathcal G}, \circ)$-stable matrix. Then obviously $\sigma({\mathbf A}) = \sigma({\mathbf L}{\mathbf A}) \subset {\mathfrak D}$.
\item[\rm 3.] Let ${\mathbf G}_0\circ{\mathbf A}$ be $({\mathfrak D},{\mathcal G},\circ)$-stable. Consider $\mathbf A$. Then
$${\mathbf G}\circ{\mathbf A} = {\mathbf G}\circ{\mathbf L}\circ{\mathbf A} = {\mathbf G}\circ({\mathbf G}_0^{-1}{\mathbf G}_0)\circ{\mathbf A} =$$
$$ = [\mbox{associativity}] = ({\mathbf G}\circ{\mathbf G}_0^{-1})\circ({\mathbf G}_0\circ{\mathbf A}) = {\mathbf G}_1\circ({\mathbf G}_0\circ{\mathbf A}), $$ where ${\mathbf G}_1 := {\mathbf G}\circ{\mathbf G}_0^{-1} \in {\mathcal G}$.
Thus $\sigma({\mathbf G}\circ{\mathbf A}) = \sigma({\mathbf G}_1\circ({\mathbf G}_0\circ{\mathbf A})) \subset {\mathfrak D}$.
\end{enumerate}
$\square$

Consider the operation of matrix addition, then the identity element is a zero matrix $\mathbf O$. For matrix multiplication, it is an identity matrix ${\mathbf I}$, for Hadamard multiplication, it is a matrix ${\mathbf E}$ with all the entries $e_{ij} = 1$.

\begin{corollary}(See, for example, \cite{LOG}, p. 79.)
\begin{enumerate}
 \item[\rm 1.] The set of (multiplicative) $D$-stable matrices is invariant under multiplication by a positive diagonal matrix $\mathbf D$.
\item[\rm 2.] Any $D$-stable matrix is stable. 
\end{enumerate}    
\end{corollary}

\subsection{Basic properties of $({\mathfrak D}, \ {\mathcal G}, \ \circ)$-stable matrices}
Here, we consider some basic matrix operations, which preserve the class of $({\mathfrak D}, \ {\mathcal G}, \ \circ)$-stable matrices for some specified stability regions ${\mathfrak D}$, matrix classes ${\mathcal G}$ and binary operations $\circ$. We mostly focus on the properties of the binary operation $\circ$.

\paragraph{Transposition} First, consider ${\mathbf A}^T$ (the transpose of ${\mathbf A}$).
\begin{theorem} Let ${\mathfrak D}\subset{\mathbb C}$ be an arbitrary (symmetric with respect to the real axes) stability region, ${\mathcal G}\subset {\mathcal M}^{n \times n}$ be an arbitrary matrix class and $\circ$ be a binary operation, satisfying the following property:
\begin{equation}\label{prop1}({\mathbf G}\circ{\mathbf A})^T = {\mathbf A}^T\circ{\mathbf G}^T \end{equation} for any matrix ${\mathbf A} \in {\mathcal M}^{n \times n}$ and any ${\mathbf G} \in {\mathcal G}$.
Then a matrix $\mathbf A$ is left $({\mathfrak D}, \ {\mathcal G}, \ \circ)$-stable if and only if ${\mathbf A}^T$ is right $({\mathfrak D}, \ {\mathcal G}^T, \ \circ)$-stable. If, in addition, $\mathcal G$ is closed with respect to the matrix transposition (i.e. ${\mathbf G} \in {\mathcal G}$ if and only if ${\mathbf G}^T \in {\mathcal G}$) and the equality
\begin{equation}\label{comm}
\sigma({\mathbf G}\circ{\mathbf A}) = \sigma({\mathbf A}\circ{\mathbf G})
\end{equation}
holds,
 then $\mathbf A$ is $({\mathfrak D}, \ {\mathcal G}, \ \circ)$-stable if and only if ${\mathbf A}^T$ is $({\mathfrak D}, \ {\mathcal G}, \ \circ)$-stable.
\end{theorem}
{\bf Proof.}  Let ${\mathbf A}$ be left $({\mathfrak D}, \ {\mathcal G}, \ \circ)$-stable. Consider ${\mathbf A}^T$. Take arbitrary $\widetilde{\mathbf G}$ from the class ${\mathcal G}^T$. Applying Property \eqref{prop1} to $\widetilde{{\mathbf G}}\circ{\mathbf A}^T$, we obtain
$${\mathbf A}^T \circ \widetilde{{\mathbf G}} = (\widetilde{{\mathbf G}}^T\circ{\mathbf A})^T, $$
which implies
$$ \sigma({\mathbf A}^T \circ\widetilde{ {\mathbf G}}) = \sigma(\widetilde{{\mathbf G}}^T\circ{\mathbf A})^T = \sigma(\widetilde{{\mathbf G}}^T\circ{\mathbf A})$$
since the spectra of real-valued matrices are symmetric with respect to the real axes. The inclusion ${\mathbf G}^T \in {\mathcal G}$ implies $\sigma({\mathbf G}^T\circ{\mathbf A}) \subset {\mathfrak D}$.

The second part of the lemma obviously follows from Property \eqref{comm}. $\square$

Subsection 1.4 shows that the operations of matrix multiplication, matrix addition and Hadamard matrix multiplication satisfy Property \eqref{prop1} and Property \eqref{comm}. All the matrix classes, mentioned in Subsection 1.3, are closed with respect to matrix transposition. Classes ${\mathcal B}_k$ and ${\mathcal B}_k^+$ are also closed with respect to matrix transposition. Thus the $({\mathfrak D}, \ {\mathcal G}, \ \circ)$-stability classes mentioned in Subsection 1.1, except, probably, Class 15, satisfy the conditions of the above theorem. For many of them, it was pointed out before, see, for example, the following corollaries.
\begin{corollary} The following statements hold:
 $\mathbf A$ is multiplicative $D$-stable if and only if ${\mathbf A}^T$ is multiplicative $D$-stable (see \cite{JOHN2}).
 \end{corollary}
 \begin{corollary}
 $\mathbf A$ is multiplicative $D$-positive if and only if ${\mathbf A}^T$ is multiplicative $D$-positive (see \cite{BAO}).
\end{corollary}
\begin{corollary}
 $\mathbf A$ is $B_k$-stable if and only if ${\mathbf A}^T$ is $B_k$-stable (see \cite{JOHD}).
\end{corollary}

\paragraph{Inversion} Now, suppose $\circ$ be associative and invertible and consider $(\circ{\mathbf A})^{-1}$ (the inverse of ${\mathbf A}$ with respect to the binary operation $\circ$).

 Assume that there are two stability regions $\mathfrak D, \widetilde{\mathfrak D} \subset {\mathbf C}$ such that
\begin{equation}\label{in}
\sigma(\mathbf A) \subset {\mathfrak D} \ \mbox{implies} \ \sigma((\circ{\mathbf A})^{-1})\subset\widetilde{\mathfrak D}.
 \end{equation}
In particular cases, we know a one-to-one mapping $\varphi: {\mathbb C} \rightarrow {\mathbb C}$ which connects the spectra of $\mathbf A$ and $(\circ{\mathbf A})^{-1}$:
 $$\sigma((\circ{\mathbf A})^{-1}) = \varphi(\sigma(\mathbf A)).$$
 The following statement holds.
    \begin{theorem}\label{inver}
Let $\circ$ be an associative and invertible matrix operation, ${\mathfrak D}$ and $\widetilde{\mathfrak D} \subset{\mathbb C}$ be two stability regions connected by Property \eqref{in}, and ${\mathcal G}\subset {\mathcal M}^{n \times n}$ be a class of invertible with respect to $\circ$ matrices. Then a matrix $\mathbf A$ is left $({\mathfrak D}, \ {\mathcal G}, \ \circ)$-stable implies $(\circ{\mathbf A})^{-1}$ is right $(\widetilde{{\mathfrak D}}, \ {\mathcal G}^{-1}, \ \circ)$-stable. If, in addition, the region $\widetilde{{\mathfrak D}}$ is also connected to the region ${\mathfrak D}$ by Property \eqref{in}, i.e.
$$ \sigma((\circ{\mathbf A})^{-1})\subset\widetilde{\mathfrak D} \ \mbox{implies} \ \sigma(\mathbf A) \subset {\mathfrak D},$$
 matrix class ${\mathcal G}\subset {\mathcal M}^{n \times n}$ is closed with respect to $\circ$-inversion and Property \eqref{comm} holds, then $\mathbf A$ is $({\mathfrak D}, \ {\mathcal G}, \ \circ)$-stable if and only if $(\circ{\mathbf A})^{-1}$ is $({\mathfrak D}, \ {\mathcal G}, \ \circ)$-stable.
\end{theorem}
{\bf Proof.} For the proof of the first part, let ${\mathbf A}$ be $({\mathfrak D}, \ {\mathcal G}, \ \circ)$-stable. Consider $(\circ{\mathbf A})^{-1}$. Taking arbitrary $(\circ\mathbf G)^{-1}$ from the class ${\mathcal G}^{-1}$, by associativity and invertibility we obtain:
$$(\circ{\mathbf A})^{-1}\circ(\circ\mathbf G)^{-1} = (\circ({\mathbf G}\circ{\mathbf A}^{-1}))^{-1}. $$
Since $\sigma({\mathbf G}\circ{\mathbf A}) \in {\mathfrak D}$, we have $\sigma((\circ{\mathbf A})^{-1}\circ(\circ\mathbf G)^{-1}) = \sigma((\circ({\mathbf G}\circ{\mathbf A}^{-1}))^{-1})\subset \widetilde{{\mathfrak D}}$.

The second part of the lemma is by applying the same reasoning and Property \eqref{comm} to $(\circ{\mathbf A})^{-1}$. $\square$

For matrix addition and matrix multiplication, we know the concrete functions $\varphi_{\cdot}(\lambda)=\frac{1}{\lambda}$ and $\varphi_+(\lambda)= -\lambda$. The regions, invariant with respect to $\varphi_{\cdot}$ are the following:
\begin{enumerate}
\item[-] the circle $\{\lambda \in {\mathbb C}: |\lambda|=1\}$;
\item[-] conic regions (both positive and negative directions of the real axes, the imaginary axes, sector regions, left and right half-planes of the complex plane, the complex plane without the imaginary axes, etc), without the origin.
\end{enumerate}

The regions, invariant with respect to $\varphi_{+}$ are:
\begin{enumerate}
\item[-] the unit disc;
\item[-] line-containing regions (the real and imaginary axes,  the complex plane without the imaginary axes, etc).
\end{enumerate}

The matrix classes ${\mathcal G}\subset {\mathcal M}^{n \times n}$ which are closed with respect to multiplicative inverse are symmetric and symmetric positive definite matrices, $\alpha$-block diagonal matrices, diagonal matrices and any fixed sign pattern class of diagonal matrices, $\alpha$-scalar matrices and vertex diagonal matrices. The classes which are closed with respect to the additive inverse are: symmetric, diagonal and vertex diagonal matrices.

Thus, the following $({\mathfrak D}, \ {\mathcal G}, \ \circ)$-stability classes satisfy the conditions of Theorem \ref{inver}: $D$-stable, $D_\alpha$-stable, $D$-positive, $D$-hyperbolic, $H$-stable, $H_\alpha$-stable matrices.

For different partial cases of the general definition of $({\mathfrak D}, \ {\mathcal G}, \ \circ)$-stability, the following observations were obtained.
\begin{corollary} $\mathbf A$ is multiplicative $D$-stable if and only if ${\mathbf A}^{-1}$ is multiplicative $D$-stable (see \cite{JOHN2}).
\end{corollary}
\begin{corollary} $\mathbf A$ is multiplicative $D$-positive if and only if ${\mathbf A}^{-1}$ is multiplicative $D$-positive (see \cite{BAO}).
\end{corollary}

\paragraph{Multiplication by a scalar} Given a finite or infinite interval $(\underline{\alpha},\overline{\alpha})$ of the real line, we refer to the properties \eqref{scaas} and \eqref{scadist} (see Section 4) connecting a binary operation $\circ$ to the operation of scalar multiplication. The following statement holds.

\begin{theorem} Let $(\underline{\alpha},\overline{\alpha}) \subseteq {\mathbb R}$ be a finite or infinite interval, ${\mathfrak D}\subset{\mathbb C}$ be a stability region satisfying $\alpha\lambda \in {\mathfrak D}$ for any $\lambda \in {\mathfrak D}$, $\alpha \in (\underline{\alpha},\overline{\alpha})$, ${\mathcal G} \subset {\mathcal M}^{n \times n}$ be an arbitrary matrix class and $\circ$ be a binary operation. If one of the following cases holds
 \begin{enumerate} \item[\rm 1.] the operation $\circ$ is connected to scalar multiplication by Property \eqref{scaas}
 \item[\rm 2.] $\circ$ is connected to scalar multiplication by Property \eqref{scadist} and the matrix class ${\mathcal G} \subset {\mathcal M}^{n \times n}$ satisfies $\frac{1}{\alpha}{\mathbf G} \in {\mathcal G}$ for any ${\mathbf G} \in {\mathcal G}$ and any $\alpha \in (\underline{\alpha},\overline{\alpha})$.
 \end{enumerate}
 then a matrix $\mathbf A$ is $({\mathfrak D}, \ {\mathcal G}, \ \circ)$-stable implies $\alpha{\mathbf A}$ is $({\mathfrak D}, \ {\mathcal G}, \ \circ)$-stable for any $\alpha \in (\underline{\alpha},\overline{\alpha})$.
\end{theorem}
{\bf Proof.}\begin{enumerate}
\item[\rm 1.] Since ${\mathbf G}\circ(\alpha{\mathbf A}) = \alpha({\mathbf G}\circ{\mathbf A})$ and $\alpha{\mathfrak D} \subseteq {\mathfrak D}$, we have
$$\sigma({\mathbf G}\circ(\alpha{\mathbf A})) = \alpha\sigma({\mathbf G}\circ{\mathbf A}) \subset {\mathfrak D}. $$
\item[\rm 2.] Since ${\mathbf G}\circ(\alpha{\mathbf A}) = \alpha((\frac{1}{\alpha}{\mathbf G})\circ{\mathbf A})$ and $\frac{1}{\alpha}{\mathbf G} \in {\mathcal G}$ and $\alpha{\mathfrak D} \subseteq {\mathfrak D}$, we have
$$\sigma({\mathbf G}\circ(\alpha{\mathbf A})) = \alpha\sigma((\frac{1}{\alpha}{\mathbf G})\circ{\mathbf A}) \subset {\mathfrak D}. $$
\end{enumerate} $\square$

Considering $(\underline{\alpha},\overline{\alpha}) = {\mathbb R}$ (except, probably, zero), we obtain that the domain ${\mathfrak D}$ consists of lines coming through the origin. Thus if $\mathbf A$ is a $D$-hyperbolic matrix, $\alpha\mathbf A$ is also $D$-hyperbolic for any $\alpha \in {\mathbb R}$. In its turn, considering $(\underline{\alpha},\overline{\alpha}) = (0, +\infty)$ (except, probably, zero), we obtain that the domain ${\mathfrak D}$ consists of half-lines coming from the origin. As examples, we may consider positive direction of the real axes, open and closed right (left) halfplane and so on. Thus if $\mathbf A$ is $D$-positive ($D$-stable), $\alpha\mathbf A$ is also $D$-positive (respectively, $D$-stable) for any $\alpha >0$. Considering $(\underline{\alpha},\overline{\alpha}) =(-1,1)$, we obtain, that the multiplication by $\alpha$ maps the unit disk ${\mathfrak D} = \{z \in {\mathbb C}: |z|< 1\}$ into itself. Thus if $\mathbf A$ is Schur $D$-stable, $\alpha\mathbf A$ is also Schur $D$-stable for any $\alpha \in (-1, 1)$.

\paragraph{Similarity transformations} Here, we consider the following question: given a nonsingular matrix $\mathbf S$ and a $({\mathfrak D},{\mathcal G}, \circ)$-stable matrix $\mathbf A$, when ${\mathbf S}{\mathbf A}{\mathbf S}^{-1}$ is again $({\mathfrak D},{\mathcal G}, \circ)$-stable? I.e. we describe the class of similarity transformation that preserve $({\mathfrak D},{\mathcal G}, \circ)$-stability.

\begin{theorem}
Given a class of nonsingular matrices ${\mathcal S} \subset {\mathcal M}^{n \times n}$ (closed with respect to multiplicative inversion), a stability region ${\mathfrak D}\subset{\mathbb C}$, a matrix class ${\mathcal G}\subset {\mathcal M}^{n \times n}$ and a binary matrix operation $\circ$. If one of the following cases holds
\begin{enumerate}
\item[\rm 1.] the operation $\circ$ is matrix multiplication and matrix class ${\mathcal S}$ commutes with the  matrix class ${\mathcal G}$;
\item[\rm 2.] If the operation $\circ$ is connected to matrix multiplication by Property \eqref{*dist2} and the matrix class ${\mathcal G}$ is invariant with respect to the linear transformations from ${\mathcal S}$
\end{enumerate}
then a matrix ${\mathbf S}{\mathbf A}{\mathbf S}^{-1}$ is $({\mathfrak D},{\mathcal G}, \circ)$-stable if and only if ${\mathbf A}$ is $({\mathfrak D},{\mathcal G}, \circ)$-stable.
\end{theorem}
{\bf Proof.} Case 1. Let ${\mathbf A}$ be multiplicative $({\mathfrak D},{\mathcal G})$-stable. Consider ${\mathbf S}{\mathbf A}{\mathbf S}^{-1}$. Since an arbitrary ${\mathbf G} \in {\mathcal G}$ commutes with an arbitrary $\mathbf S \in {\mathcal S}$, we have
$${\mathbf G}{\mathbf S}{\mathbf A}{\mathbf S}^{-1} = {\mathbf S}{\mathbf G}{\mathbf A}{\mathbf S}^{-1} = {\mathbf S}({\mathbf G}{\mathbf A}){\mathbf S}^{-1}. $$
Since $$\sigma({\mathbf S}({\mathbf G}{\mathbf A}){\mathbf S}^{-1}) = \sigma({\mathbf G}{\mathbf A}) \subset {\mathfrak D}$$
we obtain that ${\mathbf S}{\mathbf A}{\mathbf S}^{-1}$ is also multiplicative $({\mathfrak D},{\mathcal G})$-stable. For the inverse direction is enough to notice that if ${\mathbf B} := {\mathbf S}{\mathbf A}{\mathbf S}^{-1}$ then ${\mathbf A} = {\mathbf S}^{-1}{\mathbf B}{\mathbf S}$ and $\mathbf G$ commutes with $\mathbf S$ if and only if $\mathbf G$ commutes with ${\mathbf S}^{-1}$.

Case 2. Let ${\mathbf A}$ be $({\mathfrak D},{\mathcal G},\circ)$-stable. Consider ${\mathbf S}{\mathbf A}{\mathbf S}^{-1}$. Applying Property \eqref{*dist2}, we obtain
$${\mathbf G}\circ{\mathbf S}{\mathbf A}{\mathbf S}^{-1} = ({\mathbf S}{\mathbf S}^{-1}{\mathbf G})\circ({\mathbf S}{\mathbf A}{\mathbf S}^{-1})= {\mathbf S}({\mathbf S}^{-1}{\mathbf G}{\mathbf S}{\mathbf S}^{-1})\circ({\mathbf A}{\mathbf S}^{-1})= $$
$$= {\mathbf S}({\mathbf S}^{-1}{\mathbf G}{\mathbf S}\circ{\mathbf A}){\mathbf S}^{-1} ={\mathbf S}(\widetilde{{\mathbf G}}\circ{\mathbf A}){\mathbf S}^{-1},$$
where $\widetilde{{\mathbf G}} = {\mathbf S}^{-1}{\mathbf G}{\mathbf S} \in {\mathcal G}$. Since $$\sigma({\mathbf S}(\widetilde{{\mathbf G}}\circ{\mathbf A}){\mathbf S}^{-1}) = \sigma(\widetilde{{\mathbf G}}\circ{\mathbf A}) \subset {\mathfrak D},$$
we obtain that ${\mathbf S}{\mathbf A}{\mathbf S}^{-1}$ is $({\mathfrak D},{\mathcal G},\circ)$-stable. The proof for the inverse direction copies the same reasoning.
$\square$

Here, let us consider several matrix classes and their commutators. As it is known, the class ${\mathcal G}$ of diagonal matrices commutes with itself and with the class of permutation matrices. Thus we obtain the following statement (see, for example, \cite{BAO}, p. 68 for the case of $D$-positive matrices, \cite{AM}, p. 450, Theorem 2 and \cite{JOHN1}, p. 54, Observation (ii), for the case of $D$-stable matrices).
\begin{corollary} Let $\mathbf A$ belong to one of the following classes: $D$-stable matrices, $D$-positive matrices, Schur $D$-stable matrices or $D$-hyperbolic matrices. Then the matrices ${\mathbf D}{\mathbf A}{\mathbf D}^{-1}$, where ${\mathbf D}$ is a diagonal matrix and ${\mathbf P}{\mathbf A}{\mathbf P}^{-1}$, where ${\mathbf P}$ is a permutation matrix, also belongs to the same class.
\end{corollary}

\subsection{Generalized diagonal stability and sufficient conditions of $({\mathfrak D},{\mathcal G},\circ)$-stability}

Here, we consider a stability region ${\mathfrak D}$, defined by generalized Lyapunov equation \eqref{GenLyap} and provide the following definition.

Given a stability region ${\mathfrak D}$, defined by Equation \eqref{GenLyap} and a subclass $\mathcal P$ of the class of symmetric positive definite matrices $\mathcal H$, we call a matrix $\mathbf A$ Volterra--Lyapunov $({\mathfrak D},{\mathcal P})$-stable, if Equation \eqref{GenLyap} admits a solution in the matrix class $\mathcal P$, i.e. if there exists a matrix ${\mathbf P} \in {\mathcal P}$ such that
$$ {\mathbf W}:=\sum_{i,j = 0}^{n-1}c_{ij}({\mathbf A}^T)^i{\mathbf P}{\mathbf A}^j$$
is positive definite.

As partial cases, we mention the following matrix classes.

\begin{enumerate}
\item[\rm 1.] An $n \times n$ real matrix $\mathbf A$ is called {\it diagonally stable} if there exists a positive diagonal matrix $\mathbf D$ such that ${\mathbf D}{\mathbf A} + {\mathbf A}^T{\mathbf D}$ is positive definite. In this case, the matrix ${\mathbf D}$ is called a {\it Lyapunov scaling factor} of ${\mathbf A}$. The concept of diagonal stability arises in \cite{QR}, referring \cite{AM} as a characterization of multiplicative $D$-stability. The property of diagonal stability is studied in \cite{CROSS} as {\it Volterra--Lyapunov stability} and in \cite{LOG} as {\it dissipativity}. For other references and tittles of this property, see \cite{LOG}, p. 82. Here, the stability region is the left hand side of the complex plane, described by classical Lyapunov equation \eqref{lyap}, the matrix class ${\mathcal P}$ is the class of positive diagonal matrices.
\item[\rm 2.] An $n \times n$ real matrix $\mathbf A$ is called {\it $\alpha$-scalar diagonally stable} if there exists a positive diagonal $\alpha$-scalar matrix ${\mathbf D}_{\alpha}$ such that ${\mathbf D}_{\alpha}{\mathbf A} + {\mathbf A}^T{\mathbf D}_{\alpha}$ is positive definite. This matrix class was introduced in \cite{HERM} under the name of Lyapunov $\alpha$-scalar stability and then studied in \cite{GUH1}, \cite{WAN}. Here, the stability region is again the left hand side of the complex plane, the matrix class ${\mathcal P}$ is the class of $\alpha$-scalar positive diagonal matrices.
\item[\rm 3.] An $n \times n$ real matrix $\mathbf A$ is called {\it $\alpha$-block diagonally stable} if there exists a block symmetric positive definite matrix ${\mathbf H}_{\alpha}$ such that ${\mathbf H}_{\alpha}{\mathbf A} + {\mathbf A}^T{\mathbf H}_{\alpha}$ is positive definite. This matrix class was mentioned in \cite{KHAK2} in connection with the study of $D_{\alpha}$-stable matrices and studied in \cite{AB1} in connection with robust stability properties. For the applications, see also \cite{KHAL2}. Here, the stability region is again the left hand side of the complex plane, the matrix class ${\mathcal P}$ is the class of $\alpha$-block symmetric positive definite matrices.
 \item[\rm 4.]  An $n \times n$ real (not necessarily symmetric) matrix $\mathbf A$ is called {\it positive definite} if its symmetric part ${\mathbf A} + {\mathbf A}^T$ is positive definite. This matrix class was introduced in \cite{JOHN6} as a generalization of positive definiteness to non-symmetric matrices. As an equivalent characterization, it was stated that ${\mathbf A} + {\mathbf A}^T$ is positive definite if and only if $x^T{\mathbf A}x > 0$ for every nonzero vector $x \in {\mathbb R}^n$. For such matrices, the term $L$-stability is also used (see \cite{GUH1}). Here, the stability region is again the left hand side of the complex plane, the matrix class ${\mathcal P}$ consists of the only one identity matrix $\mathbf I$.
 \item[\rm 5.] An $n \times n$ real matrix $\mathbf A$ is called {\it Schur diagonally stable} if there exists a positive diagonal matrix ${\mathbf D}$ such that ${\mathbf D} - {\mathbf A}^T{\mathbf D}{\mathbf A}$ is positive definite. This definition was given in \cite{BHK}. Here, the stability region is the unit disk, defined by the Stein equation \eqref{ste}, the matrix class ${\mathcal P}$ is the class of positive diagonal matrices.
\end{enumerate}

In connection with the definition of generalized Volterra--Lyapunov stability, the following crucial question arises: given a Lyapunov stability region $\mathfrak D$, two matrix classes ${\mathcal P} \subset {\mathcal H}$ and ${\mathcal G} \subset {\mathcal M}^{n \times n}$, and a binary operation $\circ$ on ${\mathcal M}^{n \times n}$, how the class of Volterra--Lyapunov $({\mathfrak D},{\mathcal P})$-stable matrices is connected to the class of $({\mathfrak D}, \ {\mathcal G}, \ \circ)$-stable matrices?

Here, we consider the following most simple cases.

\begin{theorem} Let ${\mathfrak D} = \{\lambda \in {\mathbb C}: {\rm Re}(\lambda)> 0\}$, ${\mathcal P}, {\mathcal G} \subseteq {\mathcal H}$ be two commuting subclasses of symmetric positive definite matrices, and $\circ$ be matrix multiplication or matrix addition. Then an $n \times n$ matrix ${\mathbf A}$ is both multiplicative and additive $({\mathfrak D}, {\mathcal G})$-stable if there exist a matrix ${\mathbf P}\in {\mathcal P}$ such that
\begin{equation}\label{lyap1} {\mathbf W} : ={\mathbf P}{\mathbf A} + {\mathbf A}^T{\mathbf P} \end{equation}
       is positive definite.
\end{theorem}
{\bf Proof.} First, let us consider the case of matrix multiplication. Let ${\mathbf W} : ={\mathbf P}{\mathbf A} + {\mathbf A}^T{\mathbf P}$ be symmetric positive definite for some ${\mathbf P}\in {\mathcal P}$. Then, multiplying Equality \eqref{lyap1} from the both sides on arbitrary ${\mathbf G} \in {\mathcal G}$, we obtain $$ {\mathbf G}{\mathbf W}{\mathbf G} : ={\mathbf G}{\mathbf P}{\mathbf A}{\mathbf G} + {\mathbf G}{\mathbf A}^T{\mathbf P}{\mathbf G}$$
     $$ {\mathbf G}{\mathbf W}{\mathbf G} : =({\mathbf G}{\mathbf P})({\mathbf A}{\mathbf G}) + ({\mathbf A}{\mathbf G})^T({\mathbf G}{\mathbf P})$$
     From properties of positive definite matrices (see Subsection 1.3) we obtain that ${\mathbf G}{\mathbf W}{\mathbf G}$ and ${\mathbf G}{\mathbf P}$ are also symmetric positive definite. Thus ${\mathbf A}{\mathbf G}$ is stable by Lyapunov theorem.

Now let $\circ$ be the operation of matrix addition. Again, let ${\mathbf W} : ={\mathbf P}{\mathbf A} + {\mathbf A}^T{\mathbf P}$ be symmetric positive definite for some ${\mathbf P}\in {\mathcal P}$. Consider $\widetilde{{\mathbf W}} : ={\mathbf P}({\mathbf A}+{\mathbf G}) + ({\mathbf A} + {\mathbf G})^T{\mathbf P}$. Then
$$\widetilde{{\mathbf W}} : ={\mathbf P}{\mathbf A}+{\mathbf P}{\mathbf G} + {\mathbf A}^T{\mathbf P} + {\mathbf G}{\mathbf P}=$$
$${\mathbf W} + ({\mathbf P}{\mathbf G}+{\mathbf G}{\mathbf P}).$$
It follows from the commutativity of positive definite classes ${\mathcal P}$ and ${\mathcal G}$ that ${\mathbf P}{\mathbf G}+{\mathbf G}{\mathbf P}$ is also symmetric positive definite (see properties in Subsection 1.3). Thus $\widetilde{{\mathbf W}}$ is symmetric positive definite as a sum of positive definite matrices.
      $\square$
\begin{corollary} Diagonally stable matrices are multiplicative $D$-stable (see \cite{AM}).
\end{corollary}
\begin{corollary} Positive definite (not necessarily symmetric) matrices are $H$-stable (see \cite{AM}, p. 449, Theorem 1, also \cite{OSS}, p. 82).
\end{corollary}
\begin{corollary} $\alpha$-scalar diagonally stable matrices are $H_{\alpha}$-stable (see \cite{HERM}, p. 45, Theorem 4.4).
\end{corollary}
\begin{corollary} $\alpha$-block diagonally stable matrices are $D_{\alpha}$-stable (see \cite{KHAK2}, also \cite{AB1}).
\end{corollary}
      For the class of positive diagonal matrices, the existence of positive diagonal solution of Lyapunov equation \eqref{lyap} is sufficient, but not necessary for $D$-stability (see, for example, \cite{JOHN1}). Johnson pointed the Lyapunov diagonal stability as the oldest sufficiend condition for $D$-stability (referring \cite{QR}).  In the case of $H$-stable matrices, Ostrowski and Schneider in \cite{OSS} proved the following sufficient condition for $H$-stability: {\it an $n \times n$ matrix ${\mathbf A}$ is $H$-stable if ${\mathbf A}$ + ${\mathbf A}^T$ is positive definite} (see \cite{OSS}, p. 82). Considering also the case of positive semidefinite and singular matrix ${\mathbf A}$ + ${\mathbf A}^T$, they provide the complete characterization of $H$-stable matrices (see \cite{OSS}, p. 82, Theorem 4, also p. 81 Theorem 3 for $H$-semistability), which also shows the proper inclusion of the class of positive definite matrices to the class of $H$-stable matrices. Analogically, Schur diagonally stable matrices form a proper subclass in the class of Schur $D$-stable matrices (see \cite{KAB}).
      
      Further results on this topic are considered in \cite{KU2}.
\section{General $({\mathfrak D},{\mathcal G},\circ)$-stability theory: applications and open problems}

\subsection{Binary matrix operations theory}
The problem of defining and studying different cases of $({\mathfrak D},{\mathcal G},\circ)$-stability mainly deals with the properties of the corresponding binary operation $\circ$. Here, we consider the following questions and problems, based on establishing elementary properties of $({\mathfrak D},{\mathcal G},\circ)$-stable matrices (see Subsection 2.2).

{\bf Problem 1}. Given a binary operation $\circ$ on the set ${\mathcal M}^{n \times n}$ of matrices with real entries, when the equality
$$\sigma({\mathbf A}\circ{\mathbf B}) = \sigma({\mathbf B}\circ{\mathbf A})$$
holds for every ${\mathbf A}, \ {\mathbf B} \in {\mathcal M}^{n \times n}$?

Here, we have the following most obvious cases.
\begin{enumerate}
\item[\rm 1.] When the operation $\circ$ is commutative, we have ${\mathbf A}\circ{\mathbf B} = {\mathbf B}\circ{\mathbf A}$ which implies $\sigma({\mathbf A}\circ{\mathbf B}) = \sigma({\mathbf B}\circ{\mathbf A})$.
\item[\rm 2.] When $\circ$ is matrix multiplication, defined on the set of nonsingular matrices. Then ${\mathbf A}{\mathbf B} = {\mathbf B}^{-1}({\mathbf B}{\mathbf A}){\mathbf B}$ implies $\sigma({\mathbf A}{\mathbf B}) = \sigma({\mathbf B}{\mathbf A})$.
\end{enumerate}

{\bf Problem 2.} Given a binary operation $\circ$ on the set ${\mathcal M}^{n \times n}$, when the equality
       $$({\mathbf A}\circ{\mathbf B})^T = {\mathbf B}^T\circ{\mathbf A}^T$$
holds for every ${\mathbf A}, \ {\mathbf B} \in {\mathcal M}^{n \times n}$?

The above equality obviously holds for matrix addition, matrix multiplication and Hadamard matrix multiplication.

{\bf Problem 3.} Let the operation $\circ$ on ${\mathcal M}^{n \times n}$ be associative and invertible. Given a matrix $\mathbf A$, we have an operation inverse $(\circ{\mathbf A})^{-1}$. Assume, we know the localization of $\sigma(\mathbf A)$ inside a stability region $\mathfrak D$: $$\sigma(\mathbf A) \subset {\mathfrak D}.$$ When we can find a stability region $\widetilde{\mathfrak D}$, dependent on $\mathfrak D$, such that $$\sigma(\circ{\mathbf A})^{-1} \subset \widetilde{\mathfrak D}?$$

More strictly, when we can find a bijective mapping $\varphi:\overline{{\mathbb C}} \rightarrow \overline{{\mathbb C}}$, which connects $\sigma(\mathbf A)$ and $\sigma(\circ{\mathbf A})^{-1}$? Such mappings are well-known for the operations of matrix addition and matrix multiplication.

{\bf Problem 4.} Given a binary operation $\circ$ on the set ${\mathcal M}^{n \times n}$, can we find a rule, connecting $\circ$ to the "usual" operations of matrix multiplication and matrix addition?

As an example, we mention {\it mixed-product property} (see \cite{TSOY}), which connects the operations of Kronecker multiplication $\otimes$ and "usual" matrix multiplication by the equality
$$({\mathbf A} \otimes {\mathbf B})({\mathbf C} \otimes {\mathbf D}) = ({\mathbf A}{\mathbf C})\otimes({\mathbf B}{\mathbf D}), $$
which holds for every ${\mathbf A}, \ {\mathbf B}, \ {\mathbf C}, \ {\mathbf D}   \in {\mathcal M}^{n \times n}$.

\subsection{Characterization of $({\mathfrak D},{\mathcal G},\circ)$-stability: open problems}
Now we consider the main problems, connected to the class of $({\mathfrak D},{\mathcal G},\circ)$-stable matrices.
\paragraph{Checking $({\mathfrak D}, \ {\mathcal G}, \ \circ)$-stability}
The following two approaches, as well as any its combinations are often used for establishing $({\mathfrak D}, \ {\mathcal G}, \ \circ)$-stability.
 \begin{enumerate}
\item[\rm 1.] Imposing some additional conditions on matrix $\mathbf A$. For some important cases, $\mathbf A$ is assumed to belong to a specific matrix class, defined by determinantal inequalities.
\item[\rm 2.] Considering some more wide or more narrow stability region ${\mathfrak D}$ or matrix class ${\mathcal G}$, to make a crossway to studying another stability type which would be easier to characterize.
\end{enumerate}

 We start with the problem of major importance: given a stability region ${\mathfrak D}$, a matrix class ${\mathcal G}$ and an operation $\circ$, how to verify if a given $n \times n$ matrix $\mathbf A$ is $({\mathfrak D},{\mathcal G},\circ)$-stable, using just a finite number of steps? Note that we deal with the classes ${\mathcal G}$, that contains an infinite number of matrices.

Let us observe the modern state of the characterization problem for the most important partial cases, listed in Subsection 1.1.

\begin{enumerate}
\item[\rm 1.] {\bf Multiplicative $D$-stable matrices}. Being raised in \cite{ENT} and studied in \cite{AM}, the problem of matrix $D$-stability characterization is of major importance, due to a lot of applications of this class to mathematical modeling in economics, biology, etc (see, for example, \cite{AM}, \cite{LOG}, \cite{QR}). $D$-stability characterization problem is mentioned in \cite{HER1} among the most important problems of matrix stability (see \cite{HER1}, pp. 162-163). However, the problem of finding simple and effective methods for establishing $D$-stability still remains open. The property of $D$-stability is not easy to verify even in the finite-dimensional case. For $n = 3$, Cain provided a complete description of real $D$-stable matrices (see \cite{CA}). For $n = 4$, a verifiable criterion of $D$-stability was proved by Kanovei and Logofet (see \cite{KL}), the case $n = 4$ is also considered in \cite{JOHN4} and recently in \cite{BURL1}. For $n = 5$ a huge-volume conditions based on the Routh--Hurwitz criterion were obtained in \cite{BURL2}. For $n > 5$, no necessary and sufficient characteristic of matrix $D$-stability is yet known. For arbitrary $n$, there are some necessary for $D$-stability conditions as well as some classes of structured matrices that are known to be $D$-stable (see \cite{JOHN1}, the review papers of \cite{HER1}, \cite{LOG}, the book \cite{KAB}). Among the approaches to $D$-stability study, besides small-dimensions analysis, we should mention qualitative stability analysis (see \cite{JOHN1}, \cite{QR}), different determinantal conditions (\cite{ENT}, \cite{JOHN1}, \cite{CARL2}, \cite{KU1}), generalizations of diagonally dominant and $M$-matrices (\cite{KIM}), study of matrix scalings (\cite{CARL2}), diagonal stability condition (\cite{AM}, \cite{JOHN1}, \cite{CROSS}), studying Hadamard products, different matrix subclasses and special forms (\cite{JOHN1}, \cite{DATTA}, \cite{CARDJ}), generalized singular value approach (\cite{JOHN5}, \cite{CFY}, \cite{LEE}). Verifiability of the existing conditions is widely discussed, new criteria and approaches appear (by spectral radius minimization in \cite{KHAL3}, structured singular value approach in \cite{CFY}, \cite{LEE}, by Kharitonov criterion in \cite{KOS}). A method of checking $D$-stability by solving a number of LMI was proposed in \cite{GEOH}.
\item[\rm 2.]{\bf Multiplicative $H$-stable matrices.} In spite of the characterization problem of multiplicative $D$-stability is still unsolved, the characterization problem of multiplicative $H$-stability has been solved by characterization (see \cite{CARL3}, \cite{CAS}).
\item[\rm 3.]{\bf $D$-positive and $D$-aperiodic matrices.} Though some necessary conditions as well as some classes of $D$-positive matrices were studied in \cite{BAO}, this characterization problem is not solved and even has not been studied in full volume.
\item[\rm 4.]{\bf Schur $D$-stable matrices.} While the study of continuous-time linear systems leads to the multiplicative $D$-stability problem the study of a discrete-time case leads to the problem of Schur $D$-stability and the corresponding characterization problem. Although the number of the corresponding literature is much less, this problem is mentioned in \cite{BHK}, \cite{PRY1}, \cite{KAB} and also is not solved generally. For qualitative approach to Schur stability, see \cite{BHK1}, for LMI methods, see \cite{OGH}.
\item[\rm 5.]{\bf $D$-hyperbolic matrices.} For this new matrix class, introduced in \cite{AB2}, though some examples and applications are considered, no systematic characterization is provided.
\item[\rm 6.]{\bf Additive $D$-stable matrices.} This matrix class is widely studied by the same methods, used for multiplicative $D$-stablity (see \cite{CROSS}, \cite{SADR}, \cite{LOG}, \cite{KAB}, \cite{KOS}, \cite{RS1}). Attempts to characterize additive $D$-stability are also not fully successful yet.
\item[\rm 7.]{\bf Multiplicative and additive $H(\alpha)$-stable matrices.} Lying "between" $H$-stable and $D$-stable matrices, this class is not characterized yet. However, for some special partitions $\alpha$, the full characterization may be provided.
\item[\rm 8.]{\bf $D(\alpha)$-stable matrices.} The above is true also for this class, which lies "between" stable and $D$-stable matrices.
\item[\rm 9.]{\bf Hadamard $H$-stable matrices.} Though Hadamard products are used to characterize (multiplicative) $D$-stability (see \cite{JOHN3}, \cite{JOHN1}) and Lyapunov diagonal stability (see \cite{KR}, \cite{GUH}), the study of Hadamard $D$-stability is a matter of further development.
\item[\rm 10.]{\bf $B_k$-stable and $B_k$-nonsingular matrices.} The characterization of these matrix classes is also an open problem. For some study, see \cite{DJD}.
\end{enumerate}

Together with the main problem, we should mention the following connected subproblems.

\paragraph{Describing new classes of $({\mathfrak D},{\mathcal G},\circ)$-stable matrices} Such classes are supposed to be characterized by some collection of easy-to-verify conditions. To obtain a new class description for the most general case, we are particularly interested in some easy-to-verify conditions of ${\mathfrak D}$-stability (for a given stability region ${\mathfrak D}$). For the exception of some well-known partial cases, this is a hard problem as is. For special cases of stability regions ${\mathfrak D}$, such as the left (right)-hand side of the complex plane, unit disc and real axes, a number of such conditions is obtained and used as a base of various $({\mathfrak D},{\mathcal G},\circ)$-stability criteria. For some classes of multiplicative $D$-stable matrices, see, for example, \cite{JOHN1}, \cite{DATTA}, for Schur $D$-stable \cite{BHK} and \cite{PRY1}, for $D$-positive \cite{BAO}, for additive $D$-stable \cite{GEA}.

\paragraph{Proving $({\mathfrak D},{\mathcal G},\circ)$-stability of a given matrix class} Using general results (even if they are known) usually requires a huge amount of computations. That is why finding sufficient conditions is particularly useful. The matrices we study arise in analyzing specific mathematical models, thus they are likely to have some specific properties (e.g. symmetric positive definite, oscillatory, stochastic, $M$-matrices). The problem of proving $({\mathfrak D},{\mathcal G},\circ)$-stability of a naturally arisen matrix class characterized by its determinantal properties leads to a variety of unsolved matrix problems connected to the problems of stability of dynamical systems. We can express them as embedding relations between the class of stable matrices and other matrix classes. The most important are the question of the stability of $P^2$-matrices, asked by Hershkowitz and Johnson in \cite{HERJ1} and the question of the stability of strictly GKK $\tau$-matrices by Holtz and Schneider (see \cite{HOLS}).

\subsection{Further development of $({\mathfrak D},{\mathcal G},\circ)$-stability theory}
By analogy with already highly developed theory for partial cases (multiplicative and additive $D$-stability, Schur $D$-stability), here we provide some concept closely related to $({\mathfrak D},{\mathcal G},\circ)$-stability with the description of related problems.
\paragraph{Total $({\mathfrak D},{\mathcal G},\circ)$-stability} Here, we recall the following definition (see \cite{KAB}, p. 35). A property of an $n \times n$ matrix ${\mathbf A}$ is called {\it hereditary} if every principal submatrix of ${\mathbf A}$ shares it. The property of $({\mathfrak D},{\mathcal G},\circ)$-stability is not hereditary even in the classical case of multiplicative $D$-stability (see \cite{LOG}). Thus we introduce the following class. Given a stability region ${\mathfrak D}$, a matrix class ${\mathcal G} \subset {\mathcal M}^{k \times k}$, $k = 1, \ \ldots, \ n$ and a binary operation $\circ$ defined on ${\mathcal M}^{k \times k}$, $k = 1, \ \ldots, \ n$, a matrix $\mathbf A$ is called {\it totally $({\mathfrak D},{\mathcal G},\circ)$-stable} if it is $({\mathfrak D},{\mathcal G},\circ)$-stable and every its principal submatrix is also $({\mathfrak D},{\mathcal G},\circ)$-stable. As in the case of multiplicative $D$-stability, this matrix class may be used for studying properties of principal submatrices of $({\mathfrak D},{\mathcal G},\circ)$-stable matrices and for establishing necessary conditions for $({\mathfrak D},{\mathcal G},\circ)$-stability. Special cases of triples $({\mathfrak D},{\mathcal G},\circ)$, for which $({\mathfrak D},{\mathcal G},\circ)$-stability implies total $({\mathfrak D},{\mathcal G},\circ)$-stability are also of interest.

The class of (multiplicative) totally $D$-stable matrices was introduced in \cite{QR} (see \cite{QR}, p. 314), referring \cite{ME}, where a necessary condition for total $D$-stability was given. For the definition and study of this class see also \cite{KAB}. This class also arises in connection with further defined robust $D$-stability (see, for example, \cite{HART}, p. 205).
\paragraph{Inertia and inertia preservers} Here, we are restricted to studying specific stability regions $\mathfrak D$ with
${\rm int}({\mathfrak D}) \neq \emptyset$ and $\overline{{\mathfrak D}} \neq {\mathbb C}$. So we have three nonempty sets: ${\rm int}({\mathfrak D})$, $\partial({\mathfrak D})$ and ${\rm int}({\mathfrak D}^c) = {\mathbb C}\setminus \overline{{\mathfrak D}}$. The {\it inertia} of a square matrix ${\mathbf A}$ (with respect to a given domain ${\mathfrak D}$) is defined as a triple $(i_+({\mathbf A}), \ i_0({\mathbf A}), \ i_-({\mathbf A}))$, where $i_+({\mathbf A})$ $(i_-({\mathbf A}))$ is the number of the eigenvalues of ${\mathbf A}$ inside (respectively, outside) ${\mathfrak D}$, $i_0({\mathbf A})$ is the number of the eigenvalues on the boundary of ${\mathfrak D}$. Counting the number of eigenvalues in a given domain is also a problem of great importance in engineering. An $n \times n$ real matrix $\mathbf A$ is called {\it $({\mathfrak D},{\mathcal G},\circ)$-inertia preserving} if
     $$(i_+({\mathbf G}\circ{\mathbf A}), \ i_0({\mathbf G}\circ{\mathbf A}), \ i_-({\mathbf G}\circ{\mathbf A})) = (i_+({\mathbf G}), \ i_0({\mathbf G}), \ i_-({\mathbf G}))$$
     for every matrix ${\mathbf G} \in {\mathcal G}$.
Let us consider the partial cases.

 In the case, when ${\mathfrak D} = \{z \in {\mathbb C}: {\rm Re}(z) > 0\}$, we consider $i_+({\mathbf A})$ $(i_-({\mathbf A}))$ to be the number of the eigenvalues with positive (respectively, negative) real parts, $i_0({\mathbf A})$ to be the number of the eigenvalues with zero real parts, i.e. on the imaginary axes. The study of inertia preservers under a multiplication by a symmetric matrix $\mathbf H$ ($\mathcal G$ to be the class of symmetric matrices and $\circ$ to be matrix multiplication) was started by Sylvester and continued by Ostrowski and Schneider \cite{OSS} (see \cite{OSS}, p. 76, Theorem 1), where key results, connecting inertia and stability were presented. This is used to characterize the class of $H$-stable matrices. Classical results on this theme were obtained by Taussky \cite{TAU1}, Carlson and Schneider \cite{CAS}. An overview of this topic is presented in \cite{DATTA1}, where the inertia with respect to the unit disk is also considered. The inertia is used for the characterization of the class of $D$-stable matrices (see \cite{DATTA1}, p. 582 and references therein).

 For the generalized stability region $\mathfrak D$ and the same class of symmetric matrices $\mathcal H$, the characterization of inertia preservers is posed as an open problem in \cite{DATTA1} (p. 593, Problem 1). The tridiagonal case was considered in \cite{CARL1}, further generalization was provided in \cite{CHEN}.

\paragraph{Robustness of $({\mathfrak D},{\mathcal G},\circ)$-stability} Now we introduce one more concept of great importance in system theory. Here, we again consider a specific type of stability regions ${\mathfrak D}$, so-called {\it Kharitonov regions} (for the definitions and properties see, for example, \cite{SOF}). A matrix $\mathbf A$ is said to be {\it robustly $({\mathfrak D},{\mathcal G},\circ)$-stable} if it is $({\mathfrak D},{\mathcal G},\circ)$-stable and remains $({\mathfrak D},{\mathcal G},\circ)$-stable for sufficiently small perturbations of ${\mathbf A}$. In other words, ${\mathbf A}$ is robustly $({\mathfrak D},{\mathcal G},\circ)$-stable if ${\mathbf A}$ is $({\mathfrak D},{\mathcal G},\circ)$-stable and there exists an $\epsilon > 0$ such that for any real-valued matrix ${\mathbf \Delta}$ with $\|{\mathbf \Delta}\| < \epsilon$, the matrix ${\mathbf A} + {\mathbf \Delta}$ is $({\mathfrak D},{\mathcal G},\circ)$-stable.

Note, that in general, $({\mathfrak D},{\mathcal G},\circ)$-stability is not a robust property, even in the classical case of multiplicative $D$-stability (see \cite{AB1} for the corresponding examples). Thus discovering sufficient conditions which lead to the classes of robustly $({\mathfrak D},{\mathcal G},\circ)$-stable matrices is of great importance.

The original definition of robust multiplicative $D$-stability was given in \cite{AB1} under the name of {\it strong $D$-stability}. However, the term "strong stability" is often used in literature (see, for example, \cite{CROSS}) for additive $D$-stability, thus, to avoid confusion, we prefer the term "robust".

The analysis for robustness of the 13 sufficient conditions of $D$-stability, presented in \cite{JOHN1}, was done by Kafri \cite{KAF}. In the papers \cite{CFY}, \cite{LEE} some conditions for robust $D$-stability in terms of structured singular values are proposed. For a specific type of robust problems for $D$- and $H$-stability see \cite{BIR}, \cite{FGR}. Persistence of diagonal stability under perturbations is studied in \cite{KAB}).

Even before the notion of robust $D$-stability was given, the set of $D$-stable matrices was studied from the topological point of view (with respect to the usual topology of ${\mathcal M}^{n \times n}$) (see the results in \cite{HART}, where necessary and sufficient condition for a $D$-stable matrix to be in the interior of this set was established, also \cite{CA1}). Geometric study of the set of $D$-stable matrices and their scalings (in the case of small dimensions) was done in \cite{TOG}. In \cite{CAMM} , it was shown that the interior of the set of $D$-stable matrices coincides with the set of robustly $D$-stable matrices. Note, that the set of Lyapunov diagonally stable matrices forms a proper inclusion to the interior of the set of $D$-stable matrices (see \cite{HART}). Thus every Lyapunov diagonally stable matrix is robustly $D$-stable, but the inverse is not correct.

 The class of robust $D(\alpha)$-stable matrices was analyzed in \cite{AB1} (p. 3, Definition 3), see also \cite{AB3}.
 In the same paper \cite{AB1} robustly $D$-hyperbolic and $D(\alpha)$-hyperbolic classes are analyzed. 

\paragraph{$\mathfrak D$-stability measurement and general $\mathfrak D$-stabilization problem}
Here, we introduce the concept and state some problems that are connected to robust $\mathfrak D$-stability. We start with the following question, asked by Hershkowitz (see \cite{HER1}, p. 162).

Given a $P$-matrix $\mathbf A$, how far is it from being stable?

 He outlined two directions for giving an answer:

\begin{enumerate}
\item[-] in terms of the width of a wedge around the negative direction of the real axes, which is free from eigenvalues;
\item[-] in terms of the inertia of $\mathbf A$ (how much eigenvalues are located in the closed left-hand side of the complex plane).
\end{enumerate}
The combination of this two approaches was used in \cite{HERB}, \cite{KEL}.

Here, we state the following more general problem.

{\bf Problem 5}. Given an arbitrary stability region ${\mathfrak D} \subset \mathbb{C}$, and a matrix $\mathbf A$ from ${\mathcal M}^{n \times n}$, how far is $\mathbf A$ from being ${\mathfrak D}$-stable?

The answer may use the combination of the following approaches:
\begin{enumerate}
\item[-] description of the new stability region ${\mathfrak D}_1$ such as ${\mathfrak D} \subseteq {\mathfrak D}_1$ and $\sigma(\mathbf A) \subset {\mathfrak D}_1$;
\item[-] counting the inertia of $\mathbf A$ with respect to the stability region $\mathfrak D$.
\end{enumerate}

Another problem, mentioned in \cite{HER1} is {\it multiplicative $D$-stabilization problem} (see \cite{HER1} p. 162, then p.170): given a square real-valued matrix $\mathbf A$, can we find a diagonal matrix $\mathbf D$ such that ${\mathbf D}{\mathbf A}$ is positive stable? Simple example with a circulant matrix shows that it is not always possible. For the results on a stabilization of matrices using a diagonal matrix, we refer to \cite{BAL}, \cite{YAR}, \cite{LOC}.

In full generality, we state this problem as follows:

{\bf Problem 6}. Given a matrix $\mathbf A$ from ${\mathcal M}^{n \times n}$, an arbitrary stability region ${\mathfrak D} \subset \mathbb{C}$, a matrix class ${\mathcal G} \subset {\mathcal M}^{n \times n}$ and a binary matrix operation $\circ$, when it is possible to find a matrix ${\mathbf G}_0 \in {\mathcal G}$ such that $\sigma({\mathbf G}_0 \circ{\mathbf A}) \subset {\mathfrak D}$?

A matrix $\mathbf A$ is called {\it $({\mathfrak D},{\mathcal G},\circ)$-stabilizable} if the answer to Problem 6 is affirmative. As it follows from the definition, the class of $({\mathfrak D},{\mathcal G},\circ)$-stable matrices belongs to the class of $({\mathfrak D},{\mathcal G},\circ)$-stabilizable matrices.
\paragraph{$({\mathfrak D},{\mathcal G},\circ)$ stability measurement and $({\mathfrak D},{\mathcal G},\circ)$-stabilization problem}
Here, we ask the following more specific question.

{\bf Problem 7.} Given a ${\mathfrak D}$-stable matrix $\mathbf A$, a matrix class ${\mathcal G} \subset {\mathcal M}^{n \times n}$ and a binary matrix operation $\circ$, how far $\mathbf A$ is from being $({\mathfrak D},{\mathcal G},\circ)$-stable? The directions of giving the answer to this question are as follows.
 \begin{enumerate}
 \item[-] Describing subclasses ${\mathcal G}_1$ of the class $\mathcal G$, such that ${\mathcal G_1} \subseteq {\mathcal G}$ (or conversely) and $\sigma({\mathbf G}\circ{\mathbf A}) \subset {\mathfrak D}$ for every ${\mathbf G} \in {\mathcal G}_1$. Note, that every $({\mathfrak D},{\mathcal G},\circ)$-stabilizable matrix can be considered as $({\mathfrak D},{\mathcal G_1},\circ)$-stable for some nonempty class ${\mathcal G_1} \subseteq {\mathcal G}$.
\item[-] Describing the new stability region ${\mathfrak D}_1$ such that ${\mathfrak D} \subseteq {\mathfrak D}_1$ and $\sigma({\mathbf G}\circ{\mathbf A}) \subset {\mathfrak D}_1$ for every ${\mathbf G} \in {\mathcal G}$;
\item[-] counting the inertia of ${\mathbf G}\circ{\mathbf A}$ with respect to the stability region $\mathfrak D$ while ${\mathbf G}$ is varying along the class ${\mathcal G}$.
\end{enumerate}
We may also use the combinations of the described above approaches.

As examples of partial multiplicative $D$-stability, we mention the classes of $D(\alpha)$-stable matrices and $D_{\theta}$-stable matrices.
\paragraph{Relations between different classes of $({\mathfrak D},{\mathcal G},\circ)$-stable matrices}
Besides of relations between different $({\mathfrak D},{\mathcal G},\circ)$-stability classes, described in Section 2, based on inclusion relations between stability regions and matrix classes, relations between classes, defined by different binary operations are of interest. In general form, the problem is stated as follows.

{\bf Problem 8.} Given two triples $({\mathfrak D}_1,{\mathcal G}_1,\circ)$ and $({\mathfrak D}_2,{\mathcal G}_2,\star)$, do the corresponding classes of $({\mathfrak D},{\mathcal G},\circ)$-stable matrices intersect? For which classes of matrices $\mathbf A$ do they coincide?

The relations between matrix classes, defined in Section 1 are investigated in various ways. We do not provide any diagrams here, just refer to the following papers. The relations between Lyapunov diagonally stable, multiplicative and additive $D$-stable matrices were first studied in \cite{CROSS}. In \cite{HER1}, p. 173, the diagram showing relations between Lyapunov diagonally stable, multiplicative and additive $D$-stable matrices, is provided. For some matrix types, different stability types, namely multiplicative and additive $D$-stability classes coincide (\cite{HER1}, p. 174). For the relations between matrix classes, we refer to \cite{LOG}, where multicomponent diagrams are presented, see also \cite{CADHJ}, p. 154, Fig 1, \cite{BHK3}, \cite{CALN}. The book \cite{KAB} provides a lot of information on this topic. The relations between Hadamard $H$-stability, Lyapunov diagonal stability and $D$-stability were first considered in \cite{JOHN3}, p. 304.
\paragraph{Further development: from matrices to other objects}
Here, we briefly mention natural generalizations of $D$-stability which arises during study of nonlinear systems (see \cite{BEF} and references therein), theory of $D$-stability for polynomial matrices (see \cite{HENBS}), recent studies of multidimensional matrices (tenzors) and so on.

\subsection{Related problems of control theory}
Here, we are not going to give a deep overview of system theory problems, but just mention the most studied and the most important of them to which the defined above concept can be applied.

\paragraph{ General $\mathfrak D$-stability problem, or matrix eigenvalue clustering}
 This problem, discussed in general and partial cases (see, for example, \cite{GUT2}) is of great importance for engineering. Knowing some strategies of eigenvalue clustering allows us to establish sufficient conditions for $({\mathfrak D},{\mathcal G},\circ)$-stability as well as to describe new classes of $({\mathfrak D},{\mathcal G},\circ)$-stable matrices. As it was shown before (see Theorem \ref{Dstab}) if matrix class ${\mathcal G}$ forms a subgroup with respect to a group operation $\circ$, any $({\mathfrak D},{\mathcal G},\circ)$-stable matrix is necessarily ${\mathfrak D}$-stable. And, inversely, studying $({\mathfrak D},{\mathcal G},\circ)$-stability leads us to discovering new classes of ${\mathfrak D}$-stable matrices. The matrix eigenvalue localization problem is immediately related to the problem of root localization of the corresponding polynomial, such as stability, hyperbolicity, or, in general case, lying inside-outside a given region of a complex plane. Collecting and studying different polynomial techniques is useful in answering questions about the behaviour of the eigenvalues of structured matrices.
\paragraph{Robustness of $\mathfrak D$-stability}
In practice, studying dynamic systems, some perturbations of a system matrix may occur, and, in general, the matrix entries may be known up to some small values (for example, caused by linearization error). One of the most important system dynamics problems (see \cite{BAR}) is as follows. Given a stability region ${\mathfrak D} \subset {\mathbb C}$ and a matrix ${\mathbf A} \in {\mathcal M}^{n \times n}$, when a perturbed matrix $\widetilde{\mathbf A} = {\mathbf A} + {\mathbf \Delta}$ is ${\mathfrak D}$-stable?
A number of papers are devoted to studying this property, called {\it robust stability}, with respect to different stability regions $\mathfrak D$ (see, for example, \cite{ASCM} for EMI regions).

Sometimes, special types of perturbations are considered or some information of ${\mathbf \Delta}$ is provided, and we have the following description of the {\it uncertain matrix} $\widetilde{\mathbf A}$:
$$\widetilde{\mathbf A} = {\mathbf A} + {\mathbf U}{\mathbf \Delta}{\mathbf V}, $$
where ${\mathbf U}$, ${\mathbf V}$ are known matrices and introduced to specify the {\it structure of uncertainty}, ${\mathbf \Delta}$ is bounded by its norm. In some cases, we can easily come from studying $({\mathfrak D},{\mathcal G},\circ)$-stability to studying robust $\mathfrak D$-stability problem with some special structure of uncertainty.

{\bf Example 1.} For the case of $({\mathfrak D},{\mathcal G}, +)$-stability (i.e. the operation $\circ$ is matrix addition), we have
to check, if $\sigma({\mathbf A} + {\mathbf G}) \subset {\mathfrak D}$ for every matrix ${\mathbf G} \in {\mathcal G}$. Thus, assuming the norm of ${\mathbf G}$ to be sufficiently small, we immediately obtain robust ${\mathfrak D}$-stability problem with a specified structure of uncertainty (from the class ${\mathcal G}$).

{\bf Example 2.} Considering multiplicative or Hadamard $({\mathfrak D},{\mathcal G})$-stability and using the distributivity law, we obtain:
$${\mathbf G}\circ{\mathbf A} = ({\mathbf I} + ({\mathbf G} - {\mathbf I}))\circ{\mathbf A} = {\mathbf A} + ({\mathbf G} - {\mathbf I})\circ{\mathbf A}.$$
Thus, assuming that $\|{\mathbf G} - {\mathbf I}\|$ is sufficiently small, we obtain that every multiplicative (respectively, Hadamard) $({\mathfrak D},{\mathcal G})$-stable matrix is robustly ${\mathfrak D}$-stable with the uncertainty structure $({\mathbf G} - {\mathbf I})\circ{\mathbf A}$.

{\bf Example 3.} Considering the operation of entry-wise maximum $\oplus_m$, we obtain the class of $({\mathfrak D},{\mathcal G}, \oplus_m)$-stable matrices, that for a specific choice of ${\mathcal G}$ can be considered as an interval matrix. Using the commutativity and distributivity laws:
$${\mathbf G}\oplus_m {\mathbf A} = {\mathbf A}\oplus_m{\mathbf G} = {\mathbf A} + (-{\mathbf A}) + ({\mathbf A}\oplus_m{\mathbf G}) = $$ $$
{\mathbf A} + ({\mathbf A} - {\mathbf A})\oplus_m({\mathbf G} - {\mathbf A}) = {\mathbf A} + {\mathbf O}\oplus_m({\mathbf G} - {\mathbf A}).$$
Thus, for small values of $\|{\mathbf G} - {\mathbf A}\|$, $({\mathfrak D},{\mathcal G}, \oplus_m)$-stability problem leads to robust ${\mathfrak D}$-stability problem with the uncertainty structure ${\mathbf O}\oplus_m({\mathbf G} - {\mathbf A})$.

\paragraph{Pole assignment by output feedback} Now let us consider the following long-standing open problem in the linear system theory (see, for example, \cite{KIM2}, \cite{FUK}, \cite{OKF}).

Given a continuous system
$$\dot{x} = {\mathbf A}x + {\mathbf B}u; \qquad x,u \in {\mathbb R}^n $$
$$y = {\mathbf C}x \qquad y \in {\mathbb R}^n$$
with an output feedback law
$$ u = {\mathbf K}y. $$
 The problem formulation is as follows: for a given stability region ${\mathfrak D} \subset {\mathbb C}$ and matrices ${\mathbf A}, {\mathbf B}, {\mathbf C} \in {\mathcal M}^{n \times n}$, is it possible to find $\mathbf K$ such that ${\mathbf A} + {\mathbf B}{\mathbf K}{\mathbf C}$ is ${\mathfrak D}$-stable?

 A matrix triple ${\mathbf A}, {\mathbf B}, {\mathbf C}$ is called {\it assignable} with respect to a region ${\mathfrak D} \subset {\mathbb C}$ if there exists a matrix $\mathbf K$ such that ${\mathbf A} + {\mathbf B}{\mathbf K}{\mathbf C}$ is ${\mathfrak D}$-stable. For the conditions of assignability and detailed analysis of the problem, see, for example, \cite{WON}.

Using the technique, given above in Example 2, we can easily come from the study of multiplicative $({\mathfrak D},{\mathcal G})$-stability to the study of the assignability of a matrix pair $({\mathbf A},({\mathbf G} - {\mathbf I}))$ for any ${\mathbf G} \in {\mathcal G}$.


\begin{thebibliography}{100}

\bibitem{AB1}
E.H. Abed, {\it Strong $D$-stability,} Systems Control Lett. {\bf 7} (1986), 207-212.

\bibitem{AB2}
E.H. Abed, {\it Singularly perturbed Hopf bifurcation,} IEEE Trans. Circuits and Systems {\bf 32} (1985), 1270-1280.

\bibitem{AB3}
E.H. Abed, S.P. Boyd, {\it Perturbation bounds for structured robust stability,} Proc. 27th IEEE Conf. Dec. Contr. (1988), 1029-1031.

\bibitem{ASCM}
M. Allouche, M. Souissi, M. Chaabane, D. Mehdi, {\it Robust $D$-stability analysis of an induction motor,} in 16th Mediterranean Conference on Control and Automation, (2008), 255-260.

\bibitem{ALT}
C. Altafini, {\it Stability analysis of diagonally equipotent matrices,} Automatica {\bf 49} (2013), 2780-2785.

\bibitem{ALAL}
Yu.A. Al'pin, V.S. Al'pina, {\it Combinatorial structure of $k$-semiprimitive matrix families,} Sbornik: Mathematics {\bf 207} (2016), pp. 639--651.

\bibitem{ANB}
B. Anderson, N. Bose, E. Jury, {\it A simple test for zeros of a complex polynomial in a sector,} IEEE Transactions on Automatic Control {\bf 19} (1974), 437--438.

\bibitem{ART}
 A. Arhangel'skii, M. Tkachenko, {\it Topological groups and related structures,} Atlantis
Press, 2008.

\bibitem{AM}
K.J. Arrow, M. McManus, {\it A note on dynamical stability,} Econometrica {\bf 26} (1958), 448-454.

\bibitem{BAP}
O. Bachelier, B. Pradin, {\it $\partial\mathcal{D}$-regularity for robust matrix root clustering}, IFAC Proceedings Volumes, {\bf 36} (2003), pp. 237-242.

\bibitem{BBM}
O. Bachelier, J. Bosche, D. Mehdi, {\it On matrix root-clustering in a combination of first order regions}, IFAC Proceedings Volumes, {\bf 39} (2006), pp. 405-410.

\bibitem{BHPM}
O. Bachelier, D. Henrion, B. Pradin, D. Mehdi, {\it Robust root-clustering of a matrix in intersections or unions of regions}, SIAM J. Control Optim., {\bf 43} (2004), pp. 1078--1093.

\bibitem{BAL}
C.S. Ballantine, {\it Stabilization by a diagonal matrix}, Proc. Amer. Math. Soc., {\bf 25}
(1970), pp. 728--734.

\bibitem{BARK}
Y.S. Barkovsky, {\it Rank-one perturbations method and differential operators of oscillatory type} (in Russian),
PhD Thesis, Rostov-on-Don, 1980.

\bibitem{BAO}
Y.S. Barkovsky, T.V. Ogorodnikova, {\it On matrices with positive and simple spectra}, Izvestiya SKNC VSH Natural sciences \textbf{4} (1987), 65--70.

\bibitem{BAY}
Y.S. Barkovsky, V.I. Yudovich, {\it The momenta problem and spectral theory of the operators}, Izvestiya SKNC VSH Natural sciences \textbf{4} (1975), 49--53.

\bibitem{BAR}
B.R. Barmish, {\it New tools for robustness of linear systems}, Macmillan, New York, 1994.

\bibitem{BHK}
A. Bhaya, E. Kaszkurewicz, {\it On discrete-time diagonal and $D$-stability}, Linear Algebra Appl., {\bf 187} (1993), pp. 87--104.

\bibitem{BHK1}
A. Bhaya, E. Kaszkurewicz, {\it Qualitative stability of discrete-time systems}, Linear Algebra Appl., {\bf 117} (1989), pp. 65--71.

\bibitem{BHK2}
A. Bhaya, E. Kaszkurewicz, {\it Control perspectives on numerical algorithms and matrix problems}, SIAM, 2006.

\bibitem{BHK3}
A. Bhaya, E. Kaszkurewicz, R. Santos {\it Characterizations of classes of stable matrices}, Linear Algebra Appl., {\bf 374} (2003), pp. 159--174.

\bibitem{BHAT2}
R. Bhatia, {\it Positive definite matrices}, Princeton University Press, 2007.

\bibitem{BH}
S. Bhattacharya, H. Chapellat, L. Keel, {\it Robust control: the parametric approach}, Prentice-Hall, New Jersey, 1995.

\bibitem{BELL}
R. Bellman, {\it Introduction to matrix analysis}, McGraw Hill, New York, 2nd edition, 1970.

\bibitem{BERW}
A. Berman, R.C. Ward, {\it Classes of stable and semipositive matrices,} Linear Algebra Appl., {\bf 21} (1978), pp. 163-174.

\bibitem{BEF}
B. Besselink, H.R. Feyzmahdavian, H. Sandberg, M. Johansson, {\it $D$-stability and delay-independent stability of monotone nonlinear systems with max-separable Lyapunov functions}, IEEE Conference on Decision and Control (CDC),
(2016), pp. 3172--3177.

\bibitem{BICJ1}
T.A. Bickart, E.I. Jury, {\it Regions of polynomial root clustering,} Journal of the Franklin Institute, {\bf 304} (1977), pp. 149-160.

\bibitem{BICJ2}
T.A. Bickart, E.I. Jury, {\it The Schwarz--Christoffel transformation and polynomial root clustering,} IFAC Proceedings Volumes {\bf 11} (1978), pp. 1171-1176.

\bibitem{BICJ3}
T.A. Bickart, E.I. Jury, {\it Polynomial root clustering,} Journal of the Franklin Institute, {\bf 308} (1979), pp. 487-496.

\bibitem{BIR}
J. Bierkens, A. Ran, {\it A singular $M$-matrix perturbed by a nonnegative rank-one matrix has positive principal minors; is it $D$-stable?} Linear Algebra Appl., {\bf 457} (2014), pp. 191--208.

\bibitem{BFG}
F. Blanchini, E. Franco, G. Giordano, {\it Determining the structural properties of a class of biological models}, 51th IEEE Conference on Decision and Control (2012), pp. 5505--5510.

\bibitem{BGF}
S. Boyd, L. El Ghaoui, E. Feron, V. Balakrishnan, {\it Linear matrix inequalities in system and control theory}, SIAM, 1994.

\bibitem{BOR}
P. Borwein, T. Erdelyi, {\it Polynomials and polynomial inequalities}, Springer, 1995.

\bibitem{BURL1}
L.A. Burlakova, {\it $D$-stable 4th-order matrices,} J. Sovremennie technologii. Systemniy analiz. Modelirovanie {\bf 1 (21)} (2009), 109-116.

\bibitem{BURL2}
L.A. Burlakova, {\it Conditions of $D$-stability of the fifth-order matrices,} in: Gerdt V.P., Mayr E.W., Vorozhtsov E.V. (eds) Computer Algebra in Scientific Computing. Lecture Notes in Computer Science, Springer, Berlin, Heidelberg {\bf 5743} (2009), pp. 54-65.

\bibitem{BUT}
P. Butkovi\u{c}, {\it Max-linear systems: theory and algorithms,} Springer-Verlag London Limited 2010.

\bibitem{CA}
B. Cain, {\it Real, $3 \times 3$, $D$-stable matrices,} J. Res. Nat. Bur. Standards Sect. B, {\bf 80B} (1976), 75–77.

\bibitem{CA1}
B. Cain, {\it Inside the $D$-stable matrices,} Linear Algebra Appl., {\bf 56} (1984), 237–243.

\bibitem{CA2}
B. Cain, {\it Convergent multiples of convergent operators,} Linear Algebra Appl., {\bf 299} (1999), 171–173.

\bibitem{CADHJ}
B. Cain, L.M. DeAlba, L. Hogben, C.R. Johnson, {\it Multiplicative perturbations of stable and convergent operators,} Linear Algebra Appl., {\bf 268} (1998), pp. 151--169.

\bibitem{CALN}
B. Cain, T.D. Lenker, S.K. Narayan, P. Vermeire {\it Classes of stable complex matrices defined via the theorems of Ger\u{s}gorin and Lyapunov,} Linear and Multilinear Algebra, {\bf 56} (2008), pp. 713--724.

\bibitem{CAMM}
P.J. Campo, M. Morari, {\it Achievable closed-loop properties of systems under decentralized control: conditions involving the steady-state gain}, IEEE Transactions on Automatic Control, {\bf 39} (1994), pp. 932--943.

\bibitem{CARL1}
D. Carlson, {\it Controllability, inertia and stability for tridiagonal matrices}, Linear Algebra Appl., {\bf 56} (1984), pp. 207--220.

\bibitem{CARL2}
D. Carlson, {\it A class of positive stable matrices}, J. Res. Nat. Bur. Standards Sect. B, {\bf 78B} (1974), pp. 1--2.

\bibitem{CARL3}
D. Carlson, {\it A new criterion for $H$-stability of complex matrices}, Linear Algebra Appl., {\bf 1} (1968), pp. 59--64.

\bibitem{CARDJ}
D. Carlson, B. Datta, C. Johnson, {\it A semi-definite Lyapunov theorem and the characterization of tridiagonal $D$-stable matrices}, SIAM J. Alg. Disc. Meth., {\bf 3} (1982), pp. 293--304.

\bibitem{CAS}
D. Carlson, H. Schneider, {\it Inertia theorems for matrices: the semidefinite case,} Journal of Math. Analysis and Appl. {\bf 6} (1963), 430--446.

\bibitem{CAU}
A.L. Cauchy, {\it Calcul des indices des fonctions,} J. \'{E}cole Polytech. {\bf 15}, 176–229 (1837) ({\OE}uvres {\bf 1}(2),
416–466).

\bibitem{CHEN}
C.-T. Chen, {\it A generalization of the inertia theorem}, SIAM J. Appl. Math., {\bf 25} (1973), pp. 158--161.

\bibitem{CFY}
J. Chen, M. Fan, Ch.-Ch. Yu, {\it On $D$-stability and structured singular values}, Systems and Control letters, {\bf 24} (1995), pp. 19--24.

\bibitem{CHIGA}
M. Chilali, P. Gahinet, {\it $H_{\infty}$ design with pole placement constraints: an LMI approach},
IEEE Transactions on Automatic Control, {\bf 41} (1996), pp. 358--367.

\bibitem{CGA}
M. Chilali, P. Gahinet, P. Apkarian, {\it Robust pole placement in LMI regions},
Proceedings of the 36th Conference on Decision and Control
San Diego, USA, 1997, pp. 1291--1296.

\bibitem{COH}
A. Cohn, {\it \"{U}ber die Anzahl der Wurzeln einer algebraischen Gleichung in einem Kreise}, Mathematische Zeitschrift, {\bf 14} (1922), pp. 110-148.

\bibitem{CHO}
D. Choi, {\it Inequalities related to partial trace and block Hadamard product}, Linear and Multilinear Algebra, {\bf 66} (2018), pp. 280-284.

\bibitem{CROSS}
G.W. Cross, {\it Three types of matrix stability}, Linear Algebra Appl., {\bf 20} (1978), pp. 253--263.

\bibitem{CURT}
M.L. Curtis, {\it Matrix groups}, Springer-Verlag, 1984.

\bibitem{DATTA}
B.N. Datta, {\it Stability and $D$-stability}, Linear Algebra Appl., {\bf 21} (1978), pp. 135--141.

\bibitem{DATTA1}
B.N. Datta, {\it Stability and inertia}, Linear Algebra Appl., {\bf 302-303} (1999), pp. 563--600.

\bibitem{DEC}
R. Descartes, {\it La G\'{e}om\'{e}trie}, Leyden France: Maire, 1637. (Translation: {\it The Geometry of Ren\'{e} Descartes}. La Salle, France: Open Court, 1925.)

\bibitem{DJD}
J. Drew, C. Johnson, P. van den Driessche, {\it Strong forms of nonsingularity}, Lin. Algebra Appl., {\bf 162-164} (1992), pp. 187--204.

\bibitem{ENT}
A.C. Enthoven, K.J. Arrow, {\it A theorem on expectations and the stability of equilibrium}, Econometrica, {\bf 24} (1956), pp. 288--293.

\bibitem{FAJ}
S.M. Fallat, C.R. Johnson, {\it Sub-direct sums and positivity classes of matrices}, Linear Algebra Appl., {\bf 288} (1999), pp. 149--173.

\bibitem{FL1}
 R. Fleming, G. Grossman, T. Lenker, S. Narayan, S.-C. Ong, {\it Classes of Schur $D$-stable matrices,} Linear Algebra Appl., {\bf 306} (2000), pp. 15--24.

 \bibitem{FL2}
 R. Fleming, G. Grossman, T. Lenker, S. Narayan, S.-C. Ong, {\it On Schur $D$-stable matrices,} Linear Algebra Appl., {\bf 279} (1998), pp. 39--50.

\bibitem{FGR}
J.H. Fourie, G.J. Groenewald, D.B. Janse van Rensburg, A.C.M. Ran, {\it Rank-one perturbations of $H$-positive real matrices,} Linear Algebra Appl., {\bf 439} (2013), pp. 653--674.

\bibitem{FUK}
K. Furuta, S.B. Kim, {\it Pole assignment in a specified disk,} IEEE Transactions on Automatic Control, {\bf AC-32} (1987), pp. 423--427.

\bibitem{GANT}
 F. Gantmacher, {\it The Theory of Matrices,} Volume 1,
Volume 2. Chelsea. Publ. New York, 1990.

\bibitem{GANT2}
F. Gantmacher, {\it Applications of the Theory of Matrices,} Dover Publications, 2005.

\bibitem{GEA}
 X. Ge, M. Arcak, {\it A sufficient condition for additive $D$-stability and application to reaction-diffusion models,} Systems and Control Letters, {\bf 58} (2009), pp. 736--741.

\bibitem{GEOH}
 J. C. Geromel, M.C. de Oliveira, L. Hsu, {\it LMI characterization of structural and robust stability,} Linear Algebra Appl., {\bf 285} (1998), pp. 69--80.

\bibitem{GUH}
M. Gumus, J. Xu, {\it A new characterization of simultaneous Lyapunov diagonal stability via Hadamard products,} Linear Algebra Appl., {\bf 531} (2017), pp. 220--233.

\bibitem{GUH1}
M. Gumus, J. Xu, {\it Some new results related to $\alpha$-stability,} Linear and Multilinear Algebra, {\bf 65} (2017), pp. 325--340.

\bibitem{GUT}
S. Gutman, {\it Matrix root clustering in algebraic regions,} International Journal of Control, {\bf 39} (1984), pp. 773--778.

\bibitem{GUT2}
S. Gutman, {\it Root clustering in parameter space,} Springer-Verlag Berlin, Heidelberg, 1990.

\bibitem{GUJU}
S. Gutman, E. Jury, {\it A general theory for matrix root-clustering in subregions of the complex plane,} IEEE Transactions on Automatic control, {\bf AC-26} (1981), pp. 853--863.

\bibitem{HAD}
 K.P. Hadeler, {\it Nonlinear diffusion equations in biology,} in Proceedings of the Conference on Differential Equations, Dundee 1976, Springer Lecture Notes.

\bibitem{HART}
 D.J. Hartfiel, {\it Concerning the interior of the $D$-stable matrices,} Linear Algebra Appl., {\bf 30} (1980), pp. 201--207.

\bibitem{HENBS}
 D. Henrion, O. Bachelier, M. \v{S}ebek, {\it $D$-stability of polynomial matrices,} International Journal of Control, {\bf 74} (2001), pp. 845--856.

\bibitem{HENG}
 D. Henrion, A. Garulli (eds), {\it Positive polynomials in control,} Springer, 2005.

\bibitem{HERMI}
C. Hermite, {\it “On the number of roots of an algebraic equation between two limits,” Extract of a letter from Mr. C. Hermite of Paris to Mr. Borchardt of Berlin}, J. Reine angew. Math., {\bf 52} (1856), pp. 39–51. Translation by P.C. Parks, Int. J. Cont. {\bf 26} (1977), pp. 183–196.

\bibitem{HER1}
 D. Hershkowitz, {\it Recent directions in matrix stability,} Linear Algebra Appl., {\bf 171} (1992), pp. 161--186.

\bibitem{HERB}
 D. Hershkowitz, A. Berman, {\it Localization of the spectra of $P$- and $P_0$-matrices,} Linear Algebra Appl., {\bf 52/53} (1983), pp. 383--397.

\bibitem{HERJ1}
 D. Hershkowitz and C.R. Johnson, {\it Spectra of matrices with $P$-matrix powers,} Linear Algebra Appl., {\bf 80} (1986), pp. 159--171.

\bibitem{HERK2}
D. Hershkowitz and N. Keller, {\it Positivity of principal minors, sign symmetry and stability}, Linear Algebra Appl., {\bf 364} (2003), pp. 105--124.

\bibitem{HERM}
D. Hershkowitz, N. Mashal, {\it $P^\alpha$-matrices and Lyapunov scalar stability}, ELA. {\bf 4} (1998), 39-47.

\bibitem{HIL}
R.D. Hill, {\it Inertia theory for simultaneously triangulable complex matrices,} Linear Algebra Appl. {\bf 2} (1969), 131-142.

\bibitem{HOLS}
O. Holtz, H. Schneider, {\it Open problems on GKK $\tau$-matrices,} Linear Algebra Appl. {\bf 345} (2002), 263-267.

\bibitem{HOJ}
 R. Horn, C.R. Johnson, {\it Topics in matrix analysis,} Cambridge University
Press, 1991.

\bibitem{HMN}
 R. Horn, R. Mathias, Y. Nakamura, {\it Inequalities for unitarily invariant norms and bilinear matrix products,} Linear and Multilinear Algebra, {\bf 30} (1991), pp. 303--314.

\bibitem{HHP}
F.-H. Hsiao, J.-D. Hwang, S.-P. Pan, {\it $D$-stability analysis for discrete uncertain time-delay systems,} Appl. Math. Lett., {\bf 11} (1998), pp. 109--114.

\bibitem{HUR}
A. Hurwitz, {\it \"{U}ber die Begingungen, unter welchen eine Gleichung nur Wurzeln mit negativoen reelen Teilen besitzt,} Math. Ann., {\bf 46} (1895), pp. 273--284 (Werke, {\bf 2}, pp. 533--545).

\bibitem{JEM}
R. Jeltsch, M. Mansour (eds), {\it Stability theory,} Hurwitz Centenary Conference, Ascona, 1995. Birkh\"{a}user Verlag, 1996.

\bibitem{JOHN6}
C.R. Johnson, {\it Positive definite matrices,} The American Mathematical Monthly {\bf 77} (1970), 259-264.

\bibitem{JOHN1}
C.R. Johnson, {\it Sufficient conditions for $D$-stability,} Journal of Economic Theory {\bf 9} (1974), 53-62.

\bibitem{JOHN3}
C.R. Johnson, {\it Hadamard products of matrices,} Linear and Multilinear Algebra {\bf 1:4} (1974), 295-307.

\bibitem{JOHN4}
C.R. Johnson, {\it Second, third and fourth order $D$-stability,} J. Research Nat. Bureau Standards USA {\bf B78(1)} (1974), 11-13.

\bibitem{JOHN5}
C.R. Johnson, {\it A characterization of the nonlinearity of $D$-stability,} Journal of Mathematical Economics {\bf 2} (1975), 87-91.

\bibitem{JOHD}
C.R. Johnson, P. van den Driessche, {\it Interpolation of $D$-stability and sign stability,} Linear and Multilinear Algebra, {\bf 23} (1988), 363-368.

\bibitem{JOHN2}
C.R. Johnson, S. Narayan, {\it When the positivity of the leading principal minors implies the positivity of all principal minors of a matrix,} Linear Algebra Appl., {\bf 439} (2013), 2934-2947.

\bibitem{JONCK}
E.A. Jonckheere, {\it Algebraic and differential topology of robust stability,} Oxford University Press, 1997.

\bibitem{JU2}
 E.I. Jury, {\it Inners and stability of dynamic systems,} 2nd edition, Florida: R.E. Krieger, 1982.

\bibitem{JU3}
 E.I. Jury, {\it Stability, root clustering and inners,} IFAC Proceedings, {\bf 5} (1972), 153--159.

\bibitem{JUA}
E.I. Jury, S.M. Ahn, {\it Symmetric and innerwise matrices for the root-clustering and root-distribution of a polynomial,} Journal of the Franklin Institute {\bf 293} (1972), pp. 433-450.

\bibitem{KAF}
W. Kafri, {\it Robust $D$-stability,} Applied Math. Letters {\bf 15} (2002), 7-10.

\bibitem{KL}
G.V. Kanovei and D.O. Logofet, {\it $D$-stability of 4-by-4 matrices}, Comput. Math. Math. Phys., {\bf 38} (1998), pp. 1369--1374.

\bibitem{KAB}
E. Kaszkurewicz, A. Bhaya, {\it Matrix diagonal stability in systems and computation}, Springer, 2000.

\bibitem{KEL}
R.B. Kellogg, {\it On complex eigenvalues of $M$- and $P$-matrices,} Numer. Math. {\bf 19} (1972), 170-175.

\bibitem{KHAL2}
H.K. Khalil, {\it Asymptotic stability of nonlinear multiparameter singularly perturbed systems,} Automatica {\bf 17} (1981), pp. 797-804.

\bibitem{KHAL3}
H.K. Khalil, {\it A new test for $D$-stability,} Journal of Economic Theory {\bf 23} (1980), pp. 120-122.

\bibitem{KHAK1}
H.K. Khalil, P.V. Kokotovic, {\it Control of linear systems with multiparameter singular perturbations,} Automatica {\bf 15} (1979), pp. 197-207.

\bibitem{KHAK2}
H.K. Khalil, P.V. Kokotovic, {\it $D$-stability and multi-parameter singular perturbation,} SIAM J. Control Optim. {\bf 17} (1979), pp. 56-65.

\bibitem{KIM2}
H. Kimura, {\it Pole assignment by output feedback: a longstanding open problem,}
 Proceedings of the 33rd IEEE Conference on Decision and Control (1994), pp. 2101-2105.

\bibitem{KIM}
Y. Kimura, {\it A note on sufficient conditions for D-stability,}
Journal of Mathematical Economics, {\bf 8} (1981), pp. 113-120.

\bibitem{KOG}
J. Kogan, {\it Robust stability and convexity,}
Springer--Verlag, 1995.

\bibitem{KOS}
A. Kosov, {\it On the D-stability and additive D-stability of matrices and Svicobians,}
Journal of Applied and Industrial Mathematics, {\bf 4} (2010), pp. 200-212.

\bibitem{KR}
J. Kraaijevanger, {\it A characterization of Lyapunov diagonal stability using Hadamard products,} Linear Algebra Appl., {\bf 151} (1991), pp. 245--254.

\bibitem{KU1}
O.Y. Kushel, {\it On a criterion of $D$-stabiity for $P$-matrices}, Special Matrices, {\bf 4} (2016), pp. 181-188.

\bibitem{KU2}
O.Y. Kushel, {\it Generalized Volterra--Lyapunov stability}, in preparation.

\bibitem{KUSHN}
H.J. Kushner, {\it Stochastic stability and control}, Academic Press, New York, 1967.

\bibitem{LLK}
C.-H. Lee, T.-H. Li, F.-C. Kung, {\it $D$-stability analysis for discrete systems with a time delay,} Systems and Control Letters {\bf 19} (1992), 213-219.

\bibitem{LEE}
J. Lee, T. Edgar, {\it Real structured singular value conditions for the strong $D$-stability,} Systems and Control Letters {\bf 44} (2001), 273-277.

\bibitem{LOC}
A. Locatelli, N. Schiavoni, {\it A necessary and sufficient condition for the stabilization of a matrix and its principal submatrices,} Linear Algebra Appl. {\bf 436} (2012), 2311-2314.

\bibitem{LOG}
D.O. Logofet, {\it Stronger-than-Lyapunov notions of matrix stability, or how "flowers" help solve problems in mathematical ecology,} Linear Algebra Appl. {\bf 398} (2005), 75-100.

\bibitem{MAO}
W.-J. Mao, {\it An LMI approach to $D$-stability and $D$-stabilization of linear discrete singular systems with state delay,} Appl. Math. and Computations {\bf 218} (2011), 1694--1704.

\bibitem{MAM}
M. Marcus, H. Minc, {\it A survey of matrix theory and matrix inequalities}, Allyn and Bacon. Inc., Boston, 1964.

\bibitem{MAR}
M. Marden, {\it Geometry of polynomials}, AMS, Providence, 1966.

\bibitem{MART}
A.A. Martynyuk, {\it Stability by Liapunov's matrix function methods with applications}, Marcel Dekker, Inc., NY, 1998.

\bibitem{MELN}
V.G. Melnikov, {\it A sweeping method for matrix root clustering,} IFAC Proceedings Volumes {\bf 44} (2011), pp. 168-171.

\bibitem{MEME}
F. Mesquine, D. Mehdi, {\it Pole assignment in LMI Regions for linear constrained control systems,} in: Proceedings of the 15th Mediterranean conference on control and automation, Athens-Greece; July 27–29; 2007.

\bibitem{ME}
L. Metzler, {\it Stability of multiple markets: the Hicks conditions,} Econometrica {\bf 13} (1945), 277--292.

\bibitem{MOMOK}
T. Mori, Y. Mori, H. Kokame, {\it Common Lyapunov function approach to matrix root clustering,} System and Control Letters, {\bf 44} (2001), 73--78.

\bibitem{OKF}
R. Okabayashi, K. Furuta, {\it Arbitrary pole assignment using dynamic compensators based on linear function observers,} Proceedings of the 37th IEEE Conference on Design and Control, Tampa, Florida, USA, {\bf TA04 10:40} (1998), 1734--1739.

\bibitem{OBG}
M.C. de Oliveira, J. Bernussou, J.C. Geromel, {\it A new discrete-time robust stability condition,} Systems and Control Letters {\bf 37} (1999), 261--265.

\bibitem{OGH}
M.C. de Oliveira, J.C. Geromel, L. Hsu {\it LMI characterization of structural and robust stability: the discrete-time case,} Linear Algebra Appl., {\bf 296} (1999), 27--38.

\bibitem{OSS}
A. Ostrowski, H. Schneider, {\it Some theorems on the inertia of general matrices,} Journal of Math. Analysis and Appl. {\bf 4} (1962), 72--84.

\bibitem{PONT}
 L. Pontrjagin, {\it Topological groups,} Princeton University Press, 1946.

\bibitem{PRY1}
O. Pryporova, {\it Types of convergence of matrices,} PhD Thesis, Iowa State University, Ames, Iowa, 2009.

\bibitem{QR}
J.P. Quirk, R. Ruppert, {\it Qualitative economics and the stability of equilibrium,} Rev. Econom. Studies {\bf 32} (1965), 311-326.

\bibitem{RAS}
Q.I. Rahman, G. Schmeisser, {\it Analytic theory of polynomials,} Clarendon Press, Oxford, 2002.

\bibitem{RS1}
I.M. Romanishin, L.A. Sinitskii, {\it On the additive $D$-stability of matrices on the basis of the Kharitonov criterion,} Mathematical Notes, {\bf 72} (2002), pp. 237-240.

\bibitem{ROU1}
E.J. Routh, {\it The advanced part of a treatise on the dynamics of a system of rigid bodies}, Macmillan and Co., London, 1884, 4th ed., pp. 168--176.

\bibitem{ROU2}
E.J. Routh, {\it A treatise on the stability of a given state of motion}, Macmillan and Co., London, 1877.

\bibitem{ROU3}
E.J. Routh, {\it The Advanced Part of a Treatise on the Dynamics of a System of Rigid Bodies. Being Part II of a Treatise on the Whole Subject. With Numerous Examples}, Dover Publications, 1955.

\bibitem{ROU4}
E.J. Routh, {\it Dynamics of a system of rigid bodies}, Macmillan, 1892.

\bibitem{RUGH}
W.J. Rugh, {\it Linear system theory}, Prentice-Hall, 1996.

\bibitem{SADR}
R.A. Satnoianu, P. van den Driessche, {\it Some remarks on matrix stability with application to Turing instability}, Linear Algebra Appl. {\bf 398} (2005), 69--74.

\bibitem{SHU}
I. Schur, {\it Bemerkungen zur Theorie der beschr\"{a}nkten Bilinearformen mit unendlich vielen Ver\"{a}nderlichen}, J. Reine Angew. Math., {\bf 140} (1911), pp. 1--28.

\bibitem{SHU2}
I. Schur, {\it \"{U}ber Potenzreihen, die im Innern des Einheitskreises beschr\"{a}nkt sind}, J. Reine Angew. Math., {\bf 147} (1917), pp. 205-232.

\bibitem{SOF}
 Y.C. Soh, Y.K. Foo, {\it Kharitonov regions: it suffices to check a subset of vertex polynomials}, IEEE Transactions on Automatic Control, {\bf 36} (1991), pp. 1102--1105.

\bibitem{STE}
 P. Stein, {\it Some general theorems on iterants}, J. Res. Nat. Bur. Standards, {\bf 48} (1952), pp. 82--83.

\bibitem{STU}
 C. Sturm, {\it Analyse d'un m\'{e}moire sur la r\'{e}solution des \'{e}quations num\'{e}riques}, Bull. Sci. Math. Ferussac., {\bf II} (1829), pp. 419--422.

\bibitem{SGH}
 Y.-J. Sun, R.-S. Gau, J.-G. Hsieh, {\it Simple criteria for sector root clustering of uncertain systems with multiple time delays}, Chaos, Solutions and Fractals, {\bf 39} (2009), pp. 65--71.


\bibitem{TAU1}
O. Taussky, {\it A remark on a theorem of Lyapunov},
Journal of Math. Analysis Appl., {\bf 2} (1961), p. 105--107.

\bibitem{TAU2}
O. Taussky, {\it Matrices $C$ with $C^n \rightarrow 0$},
Journal of Algebra, {\bf 1} (1964), p. 5--10.

\bibitem{TIM}
D. Timotin, {\it Redheffer products and characteristic functions},
Journal of Mathematical Analysis Appl., {\bf 196} (1995), p. 823--840.

\bibitem{TOG}
Y. Togawa, {\it A geometric study of the $D$-stability problem},
Linear Algebra Appl., {\bf 33} (1980), p. 133--151.

\bibitem{TK}
 W. Tru\"{o}l, F.J. Kraus, {\it Robust $D$-stability in frequency domain with Kharitonov-like properties},
in: IFAC Design Methods of Control Systems, Zurich (1991), pp. 149-154.

\bibitem{TSA}
 M. Tsatsomeros, {\it Generating and detecting matrices with positive principal minors},
Asian Information-Science-Life, {\bf 1} (2002), p. 115--132.

\bibitem{TSOY}
T.-W. Ma, {\it Classical analysis on normed spaces}, World Scientific Publishing, 1995.

\bibitem{TZO}
 M. Tzoumas, {\it On sign-symmetric circulant matrices},
Applied Math. and Computations, {\bf 195} (2008), p. 604--617.

\bibitem{WAN}
 M. Wanat, {\it The $\alpha$-scalar diagonal stability of block matrices},
Linear Algebra Appl., {\bf 414} (2006), p. 304--309.

\bibitem{WAR}
E. Waring, {\it Problems}, Philos, Trans. Roy. Soc. London, {\bf 53} (1763), p. 294--299.

\bibitem{WIM}
 H.K. Wimmer, {\it Generalizations of theorems of Lyapunov and Stein},
Linear Algebra Appl., {\bf 10} (1975), p. 139--146.

\bibitem{WON}
W.M. Wonham, {\it Linear multivariable control: a geometric approach}, Third edition, Springer-Verlag, 1985.

\bibitem{YAR}
 R. Yarlagadda, {\it Stabilization of matrices},
Linear Algebra Appl., {\bf 21} (1978), p. 271--288.
\end{thebibliography}
\end{document}